\documentclass[11pt]{amsart} 

   \usepackage[a4paper, margin=100pt]{geometry}     
     \normalsize 
    \addtocontents{toc}{\setcounter{tocdepth}{1}} 
  
    \usepackage{xcolor}
    \usepackage{hyperref}
         \hypersetup{colorlinks, citecolor=blue, filecolor=blue, linkcolor=black, urlcolor=blue}
    \usepackage{soul}
    \usepackage[numbers]{natbib}    
    \usepackage{mathtools}
    \usepackage{tikz-cd}
    \usepackage{bigints}
    \usepackage{xfrac}
    \usepackage{mathrsfs}
    \usepackage{comment}
    \usepackage{verbatim}
    \usepackage{adjustbox}
    \usepackage{graphicx}
    \usepackage{import}
 
    \usepackage{marginnote}

\newcommand{\orbit}{\operatorname{Orb}}

\newcommand{\gen}[1]{\left\langle{#1}\right\rangle}

\DeclareMathOperator{\diff}{\rm{Diff}}

\DeclareMathOperator{\IFS}{\rm{IFS}}

\newcommand{\R}{\mathbb{R}}\newcommand{\N}{\mathbb{N}}

\newcommand{\T}{\mathbb{T}}

\renewcommand{\SS}{\mathbb{S}}

\newcommand{\id}{\mathrm{Id}}

\newcommand{\cB}{\mathcal{B}}
\newcommand{\cC}{\mathcal{C}}
\newcommand{\cD}{\mathcal{D}}
\newcommand{\cE}{\mathcal{E}}
\newcommand{\cF}{\mathcal{F}}
\newcommand{\cG}{\mathcal{G}}
\newcommand{\cH}{\mathcal{H}}

\newcommand{\cU}{\mathcal{U}}
\newcommand{\cV}{\mathcal{V}}
\newcommand{\cW}{\mathcal{W}}

\newcommand{\bT}{\mathbb{T}}

\numberwithin{equation}{section}

\newtheorem{theorem}{Theorem}[section] 
\newtheorem{corollary}[theorem]{Corollary}
\newtheorem{lemma}[theorem]{Lemma}
\newtheorem*{lemma*}{Lemma}
\newtheorem{proposition}[theorem]{Proposition}
\newtheorem*{proposition*}{Proposition}

\newtheorem{question}[theorem]{Question}
\newtheorem{claim}{Claim}
\newtheorem*{question*}{Question}
\newtheorem*{theorem*}{Theorem}
\newtheorem*{claim*}{Claim}

\newtheorem{theoremain}{Theorem}

\newtheorem{corollarymain}[theoremain]{Corollary}

\theoremstyle{definition}
\newtheorem{definition}[theorem]{Definition}

\newtheorem{example}[theorem]{Example}

\theoremstyle{remark}
\newtheorem{remark}[theorem]{Remark}

    \newcommand{\genplus}[1]{\langle\mathcal{#1}\rangle^{+}}
    \newcommand{\ifs}[1]{{\mathrm{IFS}(\mathcal{#1})}}
    \newcommand{\ball}[2]{B(#1,#2)}
    \newcommand{\Eucball}[1]{B_{#1}(0)}
    \newcommand{\leb}{\mathrm{Leb}}
    \newcommand{\uppercont}{\overline{\lambda}}
    \newcommand{\lowercont}{\underline{\lambda}}
    \newcommand{\upperexp}{\overline{\eta}}
    \newcommand{\lowerexp}{\underline{\eta}}
    \newcommand{\ettainemma}{\xi_1}
    \newcommand{\ozetatwo}{\kappa^2}
    \newcommand{\uzetatwo}{\kappa^{-2}}
    \newcommand{\ozeta}{\kappa^2}

    \newcommand{\oD}{\tilde{D}}
    
    \newcommand{\down}{\mathord{\downarrow}}
    \newcommand{\loc}[1]{\down_{#1}}
    \newcommand{\pd}{\hat{D}}
    \newcommand{\sldr}{\mathrm{SL}(d,\R)}
    \newcommand{\hs}[1]{\|#1\|_{2}}
    \newcommand{\matrixexp}{\mathfrak{exp}}
    \newcommand{\SOd}{\mathrm{SO}(d+1)}
    \newcommand{\sd}{\SS^d}
    \newcommand{\diag}{\mathrm{diag}}
    \newcommand{\ef}{\hat{\mathcal{D}}}
    \newcommand{\diffloc}[1]{\diff^{#1}_{\rm loc}}
    \newcommand{\pmsldr}{\mathrm{SL}^{\pm}(d,\R)}
    \newcommand{\gap}{\;:\;}

\begin{document}

    \title[Stable local ergodicity]{Stable local dynamics: \\expansion, quasi-conformality and ergodicity}
   
     \author[A. Fakhari]{Abbas Fakhari}
    \address{Department of Mathematics, Shahid Beheshti University, Tehran, 19839, Iran; \newline\indent School of Mathematics, Institute for Research in Fundamental Sciences (IPM), P.O. Box 19395-5746, Tehran, Iran}
    \email{a\_fakhari@sbu.ac.ir} 
        \author[M. Nassiri]{Meysam Nassiri}
\address{School of Mathematics, Institute for Research in Fundamental Sciences (IPM), P.O. Box 19395-5746, Tehran, Iran}
    \email{nassiri@ipm.ir}
        \author[H. Rajabzadeh]{Hesam Rajabzadeh}
    \address{School of Mathematics, Institute for Research in Fundamental Sciences (IPM), P.O. Box 19395-5746, Tehran, Iran}
    \email{rajabzadeh@ipm.ir}
    \date{}

 \thanks{
 This work is partially supported by INSF grant no.  4001845.}

    \begin{abstract}
     In this paper, we study stable ergodicity of the action of groups of diffeomorphisms on smooth manifolds. 
     Such actions are known to exist only on one-dimensional manifolds.
    The aim of this paper is to introduce a geometric method to overcome this restriction and 
     to construct
     higher dimensional examples.  In particular, we show that every closed manifold admits stably ergodic finitely generated group actions by diffeomorphisms of class $\cC^{1+\alpha}$. We also prove the stable ergodicity of certain algebraic actions,  including the natural action of a generic pair of matrices near the identity on a sphere of arbitrary dimension. These are consequences of the {\it quasi-conformal blender}, a local and stable mechanism/phenomenon introduced in this paper, which encapsulates our method for proving 
     stable local ergodicity by providing quasi-conformal orbits with fine controlled geometry. The quasi-conformal blender is developed in the context of pseudo-semigroup actions of locally defined smooth diffeomorphisms, which  allows for applications in diverse settings.
    \end{abstract}

        \maketitle

    \tableofcontents

    \section{Introduction} 
    Let $G$ be a semigroup in $\diff^{1+\alpha}(M)$,  the space of all diffeomorphisms with $\alpha$-H\"older derivative on a smooth Riemannian manifold $M$, endowed with the $\cC^1$ topology. The action of $G$ is {\it minimal} if every orbit is dense. 
    Also, the action of $G$ is {\it ergodic (w.r.t.  Leb.)} if every $G$-invariant set of positive Lebesgue measure in $M$ has full measure. Recall that a measurable set $S\subseteq  M$ is called {\it $G$-invariant}  if $f(S)\subseteq S$ up to a set of zero Lebesgue measure, for every $f\in G$. This definition of ergodicity concerns only the class of the Lebesgue measure, and the invariance of the Lebesgue measure is not assumed. A direct consequence of ergodicity is that Lebesgue almost every point has a dense orbit. We say that the action of $G$ on $M$ is {\it stably ergodic}  if there exists a finite set $\{f_1,\ldots,f_n\}\subseteq G$ such that the action of the semigroup generated by $\{\tilde{f}_1,\ldots,\tilde{f}_n\}$ is ergodic, where $\tilde{f}_i$ is   $\cC^1$ small perturbation of $f_i$ in  $\diff^{1+\alpha}(M)$. See Section \ref{sec:def}.
     
     A remarkable result regarding  ergodicity  with respect to the Lebesgue measure (in
the absence of an invariant measure in its class), proved by Katok and Herman, asserts that
every $\cC^2$ minimal circle diffeomorphism is ergodic (cf. \cite{Herman1979,Katok-Hasselblot,Navas_Book}). However, the ergodicity in this result is not stable under perturbation.  
The only known examples of stably ergodic actions by diffeomorphisms are in dimension one which come from Sullivan’s exponential expansion strategy (see for instance \cite{Shub_Sullivan_1985} and \cite[Theorem 1.6.]{Deroin_Kleptsyn_Navas_2018}):

    \begin{theorem}[Sullivan] \label{thm:sullivan} 
    Let $G$ be a group of $\cC^{1+\alpha}$ circle diffeomorphisms,  with $\alpha> 0$. Assume that for all $x \in \mathbb{S}^{1}$, there exists some $g\in G$ such that $g'(x) > 1$. If the action of $G$ is minimal,  then it is ergodic.
    \end{theorem} 

    Indeed, the action of $G$ in this theorem is stably ergodic. This important fact can be deduced from the arguments in \cite{Navas_2004SurLG}, where a proof of Theorem \ref{thm:sullivan} is given.
    A notable example of a stably ergodic group action given by Theorem \ref{thm:sullivan} is the action of a group generated by an irrational rotation and a Morse-Smale diffeomorphism.

    The aim of this paper is to introduce a mechanism for stable ergodicity, generalizing Sullivan's exponential expansion strategy to every dimension.
    In dimension one,  every action is conformal,  i.e.  balls are mapped to balls. This is a crucial fact in the study of group actions on a one-dimensional manifold and particularly in Theorem \ref{thm:sullivan} above (cf. \cite{Deroin_Kleptsyn_Navas_2018,Navas_ICM} and their references  for recent developments).   The idea of generalizing Sullivan's strategy to higher dimensions is, of course, not new. For instance, generalizations of Theorem \ref{thm:sullivan} are proved for conformal actions in higher dimensions (cf. \cite{Deroin_Kleptsyn_GAFA2007,BFMS}). 
    However, such generalizations do not provide stably ergodic actions, since  conformality 
    (or quasi-conformality) of an action is not stable in higher dimensions.

The method developed in this paper
allows us to control the geometry of images of small open balls under certain orbit-branches, in a very stable way.  In particular, it yields the following variant of Theorem \ref{thm:sullivan}. Such controlled geometry along the orbits is a crucial part of our proofs and may be of its own independent interest.

    \begin{theoremain} \label{thm:any-dim-sullivan} 
    Let $G$ be a group of $\cC^{1+\alpha}$ diffeomorphisms on a smooth closed Riemannian manifold $M$,  with $\alpha> 0$. Assume that for any $(x,v)\in T^1M$, there exists some $g\in G$ such that $m(D_xg)>1$ and  $\|\hat{D}_xg|_{v^{\perp}}\|<1$.
    If the action of $G$ is  {minimal}, then it is stably ergodic. 
    \end{theoremain}
       
    Here, $\hat{D}_x(f):=\sigma^{-1/d} D_xf$ denotes the normalized derivative, where $\sigma>0$ is the Jacobian of $f$ at $x$. Also, $m(.)$ is the co-norm of a linear map, and $v^{\perp}$ is the linear subspace orthogonal to a vector $v$. If $\dim(M)=1$, the statement of this theorem is exactly 
    that of Theorem \ref{thm:sullivan}, since $\|\hat{D}_xg|_{v^{\perp}} \|\equiv 0$ on $T^1M$. In particular, it gives another proof for the stability of ergodicity in Theorem \ref{thm:sullivan}.
    The conclusion of Theorem \ref{thm:any-dim-sullivan} is not sensitive to the choice of Riemannian metric and it suffices to ensure its assumptions for a metric on $M$. 
      
    Theorem \ref{thm:any-dim-sullivan} and its counterpart for the pseudo-groups of locally defined diffeomorphisms (Theorem \ref{thm:local-sullivan}) are obtained from a local and stable mechanism for ergodicity introduced in Theorem \ref{thm:blender}. It also allows us to prove the following. This resembles the main result of \cite{homburg_nassiri_2014} in which a similar statement is proved for minimality of semigroups action.

       \begin{theoremain}\label{thm:example} 
 Every closed Riemannian manifold $M$ admits a pair of diffeomorphisms in $\diff^{1+\alpha}(M)$, generating a stably ergodic semigroup action.
       \end{theoremain}

    This theorem, as far as the authors know, gives the first example of a stably ergodic finitely generated group action in $\diff^{1+\alpha}(M)$ on a manifold $M$ of dimension greater than one (cf. Remark \ref{rem:no-ex}). The notion of stable ergodicity in the space $\diff^{1+\alpha}(M)$ should not be mistaken with the notion of stable ergodicity within the class of volume-preserving diffeomorphisms $\diff^{1+\alpha}_{\rm vol}(M)$. In fact,  on a closed surface $S$,  we observe that the action of any cyclic subgroup of $\diff^{1+\alpha}(S)$ is not stably ergodic,  while area-preserving Anosov diffeomorphisms of $\T^2$ are stably ergodic in $\diff^{1+\alpha}_{\rm vol}(\T^2)$. 
    Such observations indicate that the number of generators in Theorem \ref{thm:example} is optimal (See Section \ref{sec:questions} for further discussion).

 On the other hand, $\cC^{1+\alpha}$ expanding endomorphisms are stably ergodic. One way to show it is by considering a dynamically defined sequence of partitions with arbitrarily small diameters, such that every element of each partition is eventually mapped to a ball of uniform size. An alternative approach in this setting is based on the functional analytic method to show the existence of a unique ergodic invariant measure in the Lebesgue class (cf.  \cite{Viana_Oliveira,Ruelle,Krzyzewski,Sacksteder_1974}). While this approach has been extended to various settings including maps with non-uniform expansions and singularities (cf. \cite{Alves_Bonatti_Viana_2000,Varandas-Viana}), it requires further developments in the setting of (pseudo) group actions.
 
 We should also mention the striking results on ergodic theory of groups of surface diffeomorphisms in \cite{Brown_Hertz_2017} which establish a classification of stationary measures for smooth group actions in dimension 2 and yield remarkable examples of stably ergodic group actions in the space of area-preserving surface diffeomorphisms  \cite{Liu,Chung_2020}. On the other hand, the stable ergodicity of finitely generated dense subgroups of isometries of even-dimensional spheres within the class of sufficiently smooth volume-preserving diffeomorphisms is proved in       \cite{Dolgopyat_Krikorian_2007} and supplemented in \cite{DeWitt}. See also \cite{zhang-etds} for similar results on transitivity. Unfortunately,  the results in these remarkable papers are not enough for showing stable ergodicity beyond the conservative setting (cf. \S\ref{subsec:st-meas}).
 
    The following algebraic example is a consequence of Theorem \ref{thm:any-dim-sullivan}. 
   
    \begin{theoremain}\label{thm:sphere-2}
    Let $d\geq 1$ and $\cF$ be a finite subset of $\mathrm{SL}(d+1,\R)$. Assume that  the closure of ${\gen{\cF}}$ strictly contains $\SOd$. Then, the natural action of $\gen{\cF}$ on $\sd$ is stably ergodic in $\diff^{1+\alpha}(\SS^d)$. Moreover, it is robustly minimal in $\diff^{1}(\SS^d)$.
    \end{theoremain}
    Here, the action of $A\in \mathrm{SL}(d+1,\R)$ on $\sd$ is defined by $x\mapsto\frac{Ax}{|Ax|}$. \\
    We say that the action of semigroup $G$ on $M$ is {\it robustly minimal} in $\diff^1(M)$, if there exists a finite set $\{f_1,\ldots,f_n\}\subseteq G$ such that the action of the semigroup generated by $\{\tilde{f}_1,\ldots,\tilde{f}_n\}$ is minimal, where $\tilde{f}_i$ is   $\cC^1$ small perturbation of $f_i$ in  $\diff^{1}(M)$ (see Section \ref{sec:def}).
 
    It is known that for $d\geq 2$, every generic pair of elements near the identity generates a dense subgroup of  $\sldr$ \cite{Kuranishi}. Therefore, every family $\cF$ formed
by a generic pair close enough to the identity satisfies the assumptions of Theorem \ref{thm:sphere-2}. Hence, we get the following corollary.
   
    \begin{corollarymain}
  \label{cor:sphere-2}
  Let $d\geq 1$. Then, for every generic pair  $A,B$ in a neighbourhood of identity  in $\mathrm{SL}(d+1,\mathbb{R})$, the natural action of $\gen{A,B}$ on $\sd$ is stably ergodic in $\diff^{1+\alpha}(\SS^d)$ and  robustly minimal in $\diff^{1}(\SS^d)$.
    \end{corollarymain}
  Let us discuss two comments on this result in the case $d=1$. In \cite{Barrientos-Raibekas}, the authors provide explicit obstructions for minimality of semigroup actions generated by diffeomorphisms of $\SS^1$ near the identity. 
  Secondly, counterexamples do exist for Corollary \ref{cor:sphere-2} when the matrices are not close to the identity (cf. Example \ref{ex:non-ergodic}). 
     
\medskip
 
    The main results of this paper imply the following uniqueness result for ergodic stationary measures. See \S\ref{sec:proof-stationary} for the definition of stationary measure.

\begin{corollarymain}\label{cor:stationary}
In the setting of either Theorems \ref{thm:any-dim-sullivan}, \ref{thm:example} or \ref{thm:sphere-2}, if the acting semigroup is finitely generated, then for any probability measure  $\nu$ on $\diff^{1+\alpha}(M)$ supported on a finite generating set, there is at most one absolutely continuous ergodic  $\nu$-stationary measure. If such a measure exists, it is equivalent to the Lebesgue measure.
\end{corollarymain}

    The problem of ergodicity is more subtle for the pseudo-groups of localized dynamics, where the restriction of the maps to a given domain is considered. While the localized dynamics appear in many settings, such as the return maps or the local holonomy of a foliation,  it is well understood only in dimension one. Most basic questions in higher dimensions are open, even for the local affine actions or local homogeneous actions (cf. \cite{Boutonnet_Ioana_Salehi} for a remarkable development for certain algebraic actions).  {As an application of our main results, one can provide the first examples of foliations of codimension greater than one which are stably ergodic with respect to the transversal Lebesgue measure (cf. \cite{RajabzadehThesis} for the detailed proof).
    }

    \subsection{Quasi-conformal blenders} \label{subsec:QC-blender}
    The next theorem provides a stable and local mechanism for generating quasi-conformal orbits of pseudo-semigroups and to deduce local ergodicity. As mentioned before, it plays a fundamental role in proving our main results stated above. Moreover, it makes the proof of Theorem \ref{thm:example} constructive, as well as flexible.

    We denote $\diffloc{1+}(M):=\bigcup_{\alpha>0}\diffloc{1+\alpha}(M)$, where $\diff_{\rm loc}^{s}(M)$ is the space of all $\cC^s$ diffeomorphisms $f : U_f \to V_f$ such that $U_f$ and $V_f$ are open subset of $M$. To obtain the strongest stability results we consider $\cC^1$ topology on this space, i.e., two elements of $\diff_{\rm loc}^s(M)$ are $\cC^1$-close if their graphs are $\cC^1$-close submanifolds of $M\times M$.
    
      A diffeomorphism $f\in\diff_{\rm loc}^{1}(M)$ is called {\it expanding} if $m(D_xf)>1$ at every point $x$ in its domain of definition.
    
    Let $\pi:\cE(M)\to M$ be the fiber bundle over a Riemannian manifold $M$ of dimension $d$ defined by
 \begin{equation}\label{def:E(M)}
        \cE(M):=\big\{(x,{\bf v}) {\gap} x\in M,~ {\bf v}\in (T_xM)^d ~ \text{and} ~ \det(A_{\bf v})=1\big\},
    \end{equation}
     where for ${\bf v}:=(v_1,\ldots,v_d)$, $A_{\bf v}$ is a $d\times d$ matrix with $(i,j)$-entries equal to $\langle v_i,v_j\rangle_x$, the inner product of $v_i,v_j$ assigned by the Riemannian metric on $T_xM$.
     
    Indeed,  $\det(A_{\bf v})$ is the volume  of the parallelogram  generated by $v_1,\ldots,v_d$ in $T_xM$. Note that $\cE(M)$ is an $\pmsldr$ bundle over $M$, where $\pmsldr$ consists of all $d\times d$ matrices with determinant $\pm1$.  Indeed, we consider a metric on $\cE(M)$ inducing the following norm on its fibers,
    \[\hs{(x,{\bf v})}:=\big(\sum\limits_{i=1}^d \langle v_i,v_i\rangle_x^2 \big)^\frac{1}{2}.\]
    In particular, whenever $M$ is an open subset of $\R^d$, $\cE(M)$ is isomorphic to the trivial bundle  $M\times \pmsldr$, endowed with the Hilbert-Schmidt norm on the fibers.

    For the oriented manifolds, where there exists a global volume form $\omega$ on $M$ compatible with the orientation, one can consider an alternative fiber bundle by replacing the condition $\det(A_{\bf v})=1$ in \eqref{def:E(M)}  with 
    $\omega|_x({\bf v})=1$. This defines a $\sldr$ fiber bundle over $M$ which is a quotient of $\cE(M)$ by an involution.

    For a diffeomorphism $f\in \diff_{\rm{loc}}^{1}(M)$ with $f:U_f\to V_f$, one can naturally define a fiber map $\ef f:\pi^{-1}(U_f)\to \pi^{-1}(V_f)$ defined by 
    \[ \ef f (x,{\bf v}):=\big(f(x),\pd_xf({\bf v})\big),\]
    where $\pd_xf({\bf v}):= ( \pd_xf(v_1),\ldots, \pd_xf(v_d))$  for ${\bf v}=(v_1,\ldots,v_d)\in (T_xM)^d$.

    Also, we use the notation $f\mathord{\downarrow}_V:= f|_{V\cap f^{-1}(V)}$  for the restriction of an invertible map $f$ to the set of points in a set $V$ that are mapped to $V$. Similarly,  we denote $\cF\mathord{\downarrow}_V:= \{f\mathord{\downarrow}_{V}  {\gap}  f\in\cF\}$ for  a family of maps localized to $V$.

    \begin{definition}[Quasi-conformal blender] \label{def:qc-blender}
         Let  $V$ be a subset of $M$ and  $\cF \subseteq \diff_{\rm loc}^{1+}(M)$ be a family of expanding diffeomorphisms between open subsets of $M$. We say $(V,\cF)$ is a {\it quasi-conformal blender}, if there exists an open subset $\cW \subseteq \cE(M)$ with compact closure such that  $V=\pi(\cW)$ and 
         \begin{equation}\label{eq:covering}
        \overline{\cW}\subseteq \bigcup\limits_{f\in \cF}(\ef f)^{-1}(\cW).
    \end{equation}

    \end{definition}
The next theorem proves local ergodicity of quasi-conformal blenders.
        
    \begin{theoremain}    \label{thm:blender} 
    Let $V\subset M$ and $\mathcal{F}\subseteq \diff^{1+}_{\rm loc}(M)$. If $(V,\cF)$ is a quasi-conformal blender, then there exists a real number $r>0$ such that 
    every measurable $\cF\mathord{\downarrow}_V$-invariant  set $S$  with $\leb(S\cap V)>0$ contains a ball of radius $r$ (up to a set of zero Lebesgue measure). Moreover, this property is $\cC^1$-stable with a uniform $r>0$.
    \end{theoremain}

    The principal role of the covering condition (\ref{eq:covering}) is to guarantee the existence of quasi-conformal orbit-branches in the domain $V$.
    In dimension one, it is equivalent to $\overline{V} \subseteq \cF^{-1}(V)$. 
   
    We would like to emphasize the multiple stability in this theorem. The covering condition (\ref{eq:covering}), the domain $V$ and the radius $r$ are all stable under small perturbations in the $\cC^1$ topology
    and are independent of the family's regularity class and its corresponding norm. This resembles the idea behind the creation of blenders in partially hyperbolic dynamics. The concept of blender was introduced in the seminal work of Bonatti and D\'iaz \cite{Bonatti_Diaz_1996} as a stable and local mechanism for transitivity. In fact, a covering condition similar to (\ref{eq:covering}) plays a key role in geometric and dynamical properties of blenders defined in \cite{Bonatti_Diaz_1996} (see \cite[Chapter 6.2]{BDV}  and \cite{homburg_nassiri_2014}).
    
    Over the last few decades, the concept of blender has been generalized and used in diverse settings (cf. \cite{HHTU,Nassiri_Pujals_2012,Avila_Crovisier_Wilkinson_2017,berger} among others).   
    A fundamental problem is to develop a variant of blender to be a stable and local mechanism for ergodicity. Indeed, Theorem  \ref{thm:blender} addresses this problem for IFS and pseudo-semigroup actions. One may expect further applications of the local tool introduced in Theorem \ref{thm:blender}  in smooth ergodic theory.

    Let us say a few words about the proof of Theorem \ref{thm:blender}.
    The proof of ergodicity is based on the simplest known method, i.e., by means of expansions.  Given a set $S$ of positive measure, one has to show that its orbit has full measure. Then, one shows that the iteration of the infinitesimal neighbourhood of a Lebesgue density point of $S$  yields large open sets that are mostly contained in the orbit of $S$.
    To implement this idea for locally defined diffeomorphisms (or for the action of diffeomorphisms that each one is expanding at some regions), controlling the geometry of balls to maintain their roundness throughout the expanding iterations is necessary.
    This necessitates the development of new techniques, which form the core contributions of this paper.    
    First, we show that covering condition (\ref{eq:covering}) implies that for some $\kappa>1$, the pseudo-semigroup generated by $\cF\mathord{\downarrow}_V$ (or its perturbations) has a $\kappa$-conformal orbit-branch at every point of $V$ (Theorem \ref{thm:conformality}).
    Secondly, we obtain a good control over the geometry of small balls under the iteration of sequences of maps satisfying infinitesimal assumptions of expansion and quasi-conformality (Theorem \ref{thm:control-shape-M}). The $\cC^{1+\alpha}$ regularity of the sequence is essential in this analysis. These two steps, in combination, demonstrate how the diversity of non-conformality can result in stable and quantitatively controlled geometry for certain iterations of small balls.
    This type of control of geometry, together with the classical distortion control argument, allows us to obtain local ergodicity. 
    Further analysis leads to uniform size of radius $r$, independent of modulus of H\"older regularity of derivatives or their norms.

    \medskip
    \noindent {\bf Organization.}
    The paper is organized as follows. In Section \ref{sec:def}, we introduce some definitions, the precise setting, and notations.
    Section \ref{sec:local} is devoted to quasi-conformality, where we show that the existence of quasi-conformal orbit-branches is equivalent to condition (\ref{eq:covering}) for some fiberwise-bounded set $\cW$. Also, two methods are given to verify that condition. They will be used in the proof of Theorems \ref{thm:any-dim-sullivan} and \ref{thm:example}. 
    In Section \ref{sec:expanding}, we prove an estimate (Theorem \ref{thm:control-shape-M}) which is a crucial technical step in the proof of all main theorems. 
    Section \ref{sec:blender} contains the proof of Theorem \ref{thm:blender} and its variants.    In Section \ref{sec:global}, we prove Theorems \ref{thm:any-dim-sullivan}, \ref{thm:example}, \ref{thm:sphere-2} and Corollary \ref{cor:stationary}. In Section \ref{sec:questions}, some questions related to the main results are discussed. 

    \medskip

\noindent{\bf Acknowledgements.}
We express our gratitude to Ali Tahzibi for his interest in this project from the beginning. We are grateful to him, Federico Rodriguez Hertz, Enrique Pujals, Daniel Smania, and Stefano Luzzatto for their valuable comments and suggestions. We also thank the anonymous referees for their comments and corrections, which have improved the readability of the paper.

    \section{Preliminary definitions and notations}\label{sec:def} 
    Let $M$ be a boundaryless smooth manifold of dimension $d$ endowed with a Riemannian metric. 
    We denote by $|.|$ the norm induced by this metric on the tangent space. We also denote the measure induced from this metric by $\leb $  and call it the Lebesgue measure. Furthermore,  we denote the ball of radius $r$ with center $x\in M$ by $\ball{x}{r}$.
    For $k\in \N$ and $\alpha\in(0,1)$,  we say $f:M\to M$ is of class $\cC^{k+\alpha}$,  whenever $f$ is $\cC^k$ and its $k$-th derivative is $\alpha$-H\"{o}lder continuous. 
    
    {For a real number $s>1$, $\cC^s$ means $\cC^{k+\alpha}$, where $s=k+\alpha$ and $k$ is its integer part. For $f\in \diffloc{1+\alpha}(\R^d)$, the $\cC^{1+\alpha}$ norm of $f$, denoted by $\|f\|_{\cC^{1+\alpha}}$ is defined by 
    \[\|f\|_{\cC^{1+\alpha}}:=\|f\|_{\cC^1}+\sup \big\{\frac{\|D_xf-D_yf\|}{|x-y|^\alpha}{\gap} x,y\in \mathrm{Dom}(f),~  0<|x-y|<1\big\}.\]
    One can use local charts and define $\cC^{1+\alpha}$ norm for the smooth maps between open subsets of manifolds.
     
    We use $\diff^s(M)$ to denote the group of $\cC^s$ diffeomorphisms of $M$. Also, $\diff^{1+}(M)$ is the union of all $\diff^s(M)$ with $s>1$. Throughout the paper,  we usually consider the $\cC^1$ topology on $\diff^s(M)$.
     Denote by $\diff_{\rm loc}^s(M)$ the space of all $\cC^s$ diffeomorphisms $f : \mathrm{Dom}(f) \to \mathrm{Im}(f)$ such that $\mathrm{Dom}(f)$, the domain of $f$ and $\mathrm{Im}(f)$, the image of $f$ are open subset of $M$. Two elements of $\diff_{\rm loc}^s(M)$ are $\cC^l$-close if their graphs  are $\cC^l$-close submanifolds of $M\times M$, for $l\leq s$. Similarly, $\diffloc{1+}(M)$ is the union of of all $\diffloc{s}(M)$ with $s>1$.}
     
    For $x\in M$,  $U\subseteq M$ and families of maps $\mathcal{F},\mathcal{G}\subseteq \diff_{\rm loc}^1(M)$,  we denote $\cF(U):=\bigcup_{f\in\cF}f(U)$,  $\cF(x):=\cF(\{x\})$,  and 
        \[\mathcal{F}\circ \mathcal{G}:=\{f\circ g {\gap}  f\in \mathcal{F},g\in \mathcal{G}\}.\]  
    Also, put $\mathcal{F}^0=\{\mathrm{Id}\}$ and for $k\in \mathbb{N}$,  denote $\mathcal{F}^k:=\mathcal{F}^{k-1}\circ \mathcal{F}$. We use $\gen{\mathcal{F}}^+$ (resp. $\gen{\mathcal{F}}$) for the pseudo-semigroup (resp. the pseudo-group) generated by $\cF$. By $\IFS{(\cF)}$,  we mean \textit{the iterated function system generated by $\cF$},  that is the action of $\gen{\cF}^+$ on $M$. 

    Given a finite family $\mathcal{F}=\{f_1,\ldots,  f_k\}\subseteq \diff_{\rm loc}^{1}(M)$, we say that the family $\tilde{\cF}\subseteq \diff_{\rm loc}^{1}(M)$ is $\epsilon$-close to $\cF$ in the $\cC^1$ topology, if $\tilde{\cF}=\{\tilde{f}_1,\ldots,\tilde{f}_k\}$ such that  $f_i,\tilde{f}_i$ are $\epsilon$-close in the $\cC^1$ topology,  for any $i=1,\ldots,k$.

    We say a property P  holds {\it $\cC^1$-stably} for $\cF$ in $\diff^{1+\alpha}(M)$,  if  P  holds   for every $\tilde{\cF}\subseteq \diff^{1+\alpha}(M)$ which is $\epsilon$-close to $\cF$ in the $\cC^1$ topology for some $\epsilon>0$.
    Also, we say a property P  holds {\it $\cC^1$-robustly} for $\cF$, if  P  holds  for every $\tilde{\cF}\subseteq \diff^{1}(M)$ which is $\epsilon$-close to $\cF$ in the $\cC^1$ topology for some $\epsilon>0$. 
    Clearly, by the definition, $\cC^1$-robustness is stronger than $\cC^1$-stability. Similarly,  $\cC^1$-stability and $\cC^1$-robustness are defined  for $\diff^{1+}(M)$, $\diff_{\rm loc}^{1+}(M)$ and $\diff_{\rm loc}^s(M)$, $s>1$.

    \subsection{Localized dynamics} 

    As it was defined in Subsection \ref{subsec:QC-blender},
    for an open set $V\subseteq M$ and $f\in\diff_{\rm loc}(M)$ with $f:\mathrm{Dom}(f)\to \mathrm{Im}(f)$, the localization of $f$ to $V$ is denoted by $f\loc{V}:=f|_{V\cap f^{-1}(V)}$.  Clearly,  the domain and the image of $f\loc{V}$ are $\mathrm{Dom}(f\loc{V})=V\cap f^{-1}(V\cap  \mathrm{Im}(V))$, and $\mathrm{Im}(f\loc{V})= f(V\cap \mathrm{Dom}(f))\cap V$, respectively, and $f\loc{V}:\mathrm{Dom}(f\loc{V})\to \mathrm{Im}(f\loc{V})$ is a bijective map.
     
    For a family $\cF$ of invertible maps, as in Subsection \ref{subsec:QC-blender}, we denote the pseudo-semigroup (resp. pseudo-group) generated by localization of elements of $\cF$ to $V$ by $\gen{\cF\loc{V}}^+$ (resp. $\gen{\cF\loc{V}}$) and by $\IFS(\cF\loc{V})$ the action of this pseudo-semigroup. A \emph{finite {orbit}-branch of $\IFS(\cF\loc{V})$ at $x$} is a sequence $\{x_i\}_{i=0}^n$ in $V$ such that $x_0=x$ and for any $1\leq i\leq  n$, there exists $f_i\in\cF$ with $x_{i-1}\in \mathrm{Dom}(f_i\loc{V})$ and $f_i(x_{i-1})=x_i$. Infinite {orbit}-branches are defined similarly. The orbit of $\IFS(\cF\loc{V})$ at $x\in V$, denoted by $\langle \cF\loc{V}\rangle^+(x)$, is the set of all points in finite  orbit-branches at $x$. For $S\subseteq V$, define $\langle \cF\loc{V}\rangle^+(S):=\cup_{x\in S}\langle \cF\loc{V}\rangle^+(x)$.

    In this paper,  we deal with two basic dynamical concepts, namely, minimality and ergodicity. For $\cF\subset \diff^1(M)$,  $\IFS(\cF\loc{V})$ is called \emph{minimal} if for any $x\in V$, $\langle \cF\loc{V}\rangle^+(x)$ is dense in $V$.  
    Throughout the paper, we consider the Lebesgue measure as the reference measure on the manifolds. A measurable set $S\subseteq V$ is called {\emph{$\cF\loc{V}$-invariant}}, if $\langle \cF\loc{V}\rangle^+(S)\subseteq S$ up to a set of  zero Lebesgue measure. Moreover, $\IFS(\cF\loc{V})$ is called {\emph{ergodic}}, if there is no  measurable $\cF\loc{V}$-invariant  set $S$ with $0<\leb(S\cap V)<\leb(V)$. 
    Note the definition of ergodicity  could be extended to actions on any measure space. However, for our specific requirements, we restrict the definitions to the ergodicity of smooth actions with respect to the Lebesgue measure.

    \subsection{Expanding and contracting maps}
    For a linear map $D$,  we denote its operator norm by $\|D \|$,  and  its co-norm by $m(D):=\inf\{|D(v)| {\gap}  |v|=1\}$. Clearly, if $D$ is a $d\times d $ matrix,
    \begin{equation}\label{eq:ineq-norm-conorm}
        m(D)\leq |\det D|^{\frac{1}{d}}\leq \|D\|,
    \end{equation}

    and if $D$ is invertible, then $m(D)=\|D^{-1}\|^{-1}$. This implies that for invertible $d\times d$ matrices $D_1,D_2,$ 
    \begin{equation}\label{eq:diff-conorm}
        \big\vert m(D_1)-m(D_2)\big\vert\leq \Big\vert \|D_1\|-\|D_2\|\Big\vert \leq \|D_1-D_2\|.
    \end{equation}

    A diffeomorphism $f\in \diffloc{1}(M)$ is called expanding, if there exists $\eta>1$ such that $m(D_xf)>\eta$ for every  $x\in {\rm Dom}(f)$. Similarly, it is called contracting, if there exists $\lambda\in (0,1)$ such that $\|D_xf\|<\lambda$ for every $x\in {\rm Dom}(f)$.
    Clearly, the expanding and contracting properties are $\cC^1$-robust.

    \begin{definition} 
    \label{def:expanding-seq} \label{def:contracting-seq} 
    For $\eta >1$ and $N\in \N\cup \{\infty\}$, we say a sequence $\{f_i\}_{i=1}^N$ in $\diff^1_{\rm loc}(M)$ is \emph{$\eta$-expanding (resp. \emph{$\eta^{-1}$-contracting})} at $x_0\in M$, if for any integer $i\in [1,N]$, $x_{i-1}\in \mathrm{Dom}(f_i)$ and $m(D_{x_{i-1}}f_i)>\eta$ (resp. $\|D_{x_{i-1}}f_i\|<\eta^{-1}$), where $x_i=f_i \circ \cdots \circ   f_1 (x_{0})$. Furthermore, the sequence is \emph{expanding at $x_0$} if it is $\eta$-expanding for some $\eta>1$.
    Similarly, the sequence is \emph{contracting at $x_0$} if it is $\lambda$-contracting for some $\lambda\in (0,1)$.
    \end{definition}

    \subsection{Quasi-conformality}\label{subsec:qc-def} For a real number $\kappa\geq 1$, a matrix $D\in \mathrm{GL}(d,\R)$ is \emph{$\kappa$-conformal}, if ${\|D\|}/{m(D)}=\|D\|.\|D^{-1}\|\leq\kappa$. A sequence $\{D_i\}_{i=1}^N$ in $\mathrm{GL}(d,\R)$ is \emph{$\kappa$-conformal}, if for any integer $n\in [1,N]$, $D_n D_{n-1}\cdots D_1$ is $\kappa$-conformal, $N$ can be finite or infinite. It follows immediately from definition that for  $D_1,D_2\in \mathrm{GL}(d,\R)$, if $D_i$ is $\kappa_i$-conformal ($i=1,2$), then $D_1^{-1}$ is $\kappa_1$-conformal and $D_1D_2$ is $\kappa_1\kappa_2$-conformal. These, in particular, imply that for a $\kappa$-conformal sequence $\{D_i\}_{i=1}^N$, all the products of the form $D_jD_{j-1}\cdots D_i$ for $1\leq i<j\leq N$ are $\kappa^2$-conformal. As derivatives of smooth maps are linear maps between tangent spaces, one can define similar notions for them. 
    \begin{definition}\label{Def:Conformal_sequence_Diffeos} 
    For $\kappa\geq 1$ and $N\in \N\cup \{\infty\}$,  we say a sequence $\{f_i\}_{i=1}^N$ in $\diff^1_{\rm loc}(M)$   is  \emph{$\kappa$-conformal at  $x\in M$},  if for any integer $n\in [1,N]$, $x\in \mathrm{Dom}(f^n)$ and the linear map $D_xf^n$ is $\kappa$-conformal, where $f^n=f_n\circ \cdots \circ f_1$. The sequence is \textit{quasi-conformal at $x$} if it is $\kappa$-conformal at $x$ for some $\kappa\geq 1$.
    \end{definition}
    
     \section{Quasi-conformal dynamics} \label{sec:local}
    This section is devoted to the notion of covering property for derivatives. We discuss its consequences in providing quasi-conformal orbit branches and also sufficient conditions to ensure it.
   
    Throughout the section, $M$ is a boundaryless Riemannian manifold of dimension $d$. We consider the fiber bundle $\pi:\cE(M)\to M$. Recall that for  ${\bf w}=(x,(v_1,\ldots,v_d) )\in \cE(M)$, $\hs{{\bf w}}$ is defined by $\hs{{\bf w}}=(\sum_{i}|v_i|^2)^\frac{1}{2}$, where $|.|$ is the norm induced by the Riemannian metric on $TM$. 
   
    Furthermore, for a linear map $T:W_0\to W_1$ between finite dimensional vector spaces endowed with inner products, we denote the Hilbert-Schmidt norm of $T$ by $\hs{T}$, defined by  \[\hs{T}:=\big(\sum_i|Te_i|^2\big)^\frac{1}{2},\]
    where $\{e_i\}_{i}$ is an orthonormal basis for     $W_0$.  In particular, for  $f\in\diff^1_{\rm loc}(M)$, and ${\bf w}=(x,(e_1,\ldots,e_d) )\in \cE(M)$, where $x\in \mathrm{Dom}(f)$ and $\{e_1,\ldots,e_d\}$ form an orthonormal basis for $T_xM$, it follows that $\hs{\pd_x f}=\hs{\ef f ({\bf w})}$.
    
    For linear isomorphisms $T:W_0\to W_1$ and $S:W_1\to W_2$ between $d$-dimensional vector spaces, we will use the following properties.
    \begin{align}
    & \|T\|\leq \hs{T}\leq \sqrt{d}\|T\|.\label{eq:norm-1}\\
    & \hs{S\circ T}\leq \|S\|.\hs{T} ~\text{and}~  \|S\circ T\|_2\leq \hs{S}.\|T\|.\label{eq:norm-2}\\
    & \text{If}~|\det T|=1,~\text{then}~\|T^{-1}\|\leq \|T\|^{d-1}.~\text{Therefore,}~T~\text{is}~ \|T\|^d\text{-conformal}. \label{eq:norm-3}\\
    & \text{If} ~ |\det T|=1~\text{and}~T~\text{is}~\kappa\text{-conformal}, ~\text{then}~ \|T\|\leq \kappa. \label{eq:norm-4}
    \end{align}
    
    \subsection{Stable quasi-conformality} 
    In this subsection, we present a criterion for the existence of quasi-conformal orbit-branches in pseudo-semigroup actions.
  
    A subset $\cW\subseteq \cE(M)$ is called {\emph{fiberwise-bounded}}, if $\sup_{{\bf w}\in \cW}\hs{{\bf w}}<\infty$.

    \begin{theorem}[Criterion for quasi-conformality] \label{thm:conformality}  Let $V\subseteq M$ be an open set and $\cF\subseteq \diff^1_{\rm loc}(M)$. Then, the following are equivalent 
    \begin{itemize}
    \item[(a)] There exists $\kappa>1$ such that for every $x\in V$, the pseudo-semigroup generated by $\cF\loc{V}$ has a $\kappa$-conformal orbit-branch at $x$.
    \item[(b)] There exists a non-empty fiberwise-bounded subset $\cW\subseteq \cE(M)$ 
    such that $\pi(\cW)=V$, and  
    \begin{equation}\label{eq:cover-W-open}
        \cW\subseteq \bigcup_{f\in \cF}(\ef f)^{-1}(\cW).
    \end{equation}
    \end{itemize}
    \end{theorem}

    \begin{proof}
    For the implication (a)$\Rightarrow$(b). For every $x\in V$, let $\{f_{x,i}\}_{i=1}^\infty$ be a sequence defining a $\kappa$-conformal orbit-branch of $\cF\loc{V}$ at $x$. Consider an arbitrary orthonormal basis $\{e_i\}_i$ for $T_xM$ and let ${\bf w}_x=(x,{\bf e})\in \pi^{-1}(x)$ be such that  ${\bf e}=(e_1,\ldots,e_d)$. Note that for any $x\in V$, $\hs{{\bf w}_x}=\sqrt{d}$. Then, the set $\cW$ defined by 
    \[\cW:=\bigcup_{x\in V}\bigcup_{i\geq 0} \ef f^i_x({\bf w}_x),\]
    satisfies the desired properties, where $f^i_x:=f_{x,i}\circ f_{x,i-1}\circ \cdots \circ f_{x,1}$. It is clear that $\pi(\cW)=V$. The covering property (\ref{eq:cover-W-open}) follows immediately from the definition of $\cW$, since every element of $\cW$ is of the following form 
    \[\ef f^i_x({\bf w}_x)=(\ef f_{x,i+1})^{-1}(\ef f_x^{i+1} ({\bf w}_x))\in (\ef f_{x,i+1})^{-1} (\cW).\]
    Finally, to show that $\cW$ is fiberwise-bounded, note that  for every $x\in V$ and $i\in \N$, $\pd_x f^i_x$ is $\kappa$-conformal and so by \eqref{eq:norm-2} and \eqref{eq:norm-4},
    \begin{equation}\label{eq:bound-norm-W}
        \hs{\ef f^i_x ({\bf w}_x)}\leq \|\pd_x f^i_x\|.\hs{{\bf w}_x}= \|\pd_x f^i_x\|\sqrt{d} \leq \kappa\sqrt{d}.
    \end{equation}
    This finishes the proof of the part (a)$\Rightarrow$(b).
    
    \medskip
    Next, in order to show (b)$\Rightarrow$(a), first denote $H:=\sup_{{\bf w}\in \cW} \hs{{\bf w}}$ and fix ${\bf w}=(x,(v_1,\ldots,v_d) )\in \cW$. We  claim that there exists a sequence $\{f_i\}_{i=1}^\infty$ in $\cF$ such that for any $n\geq 1$, $\ef f^n({\bf w})\in \cW$. The proof is by induction, assuming that $f_1,\ldots,f_n$ are defined, satisfying the properties. Then, 
    \[\ef f^n ({\bf w}) \in \cW\subseteq \bigcup_{f\in \cF}(\ef f)^{-1}(\cW).\]
    So, there is $f_{n+1}\in \cF$ with $\ef f^n({\bf w})\in (\ef f_{n+1})^{-1}(\cW)$. Consequently, by the chain rule, $\ef f^{n+1}({\bf w})\in \cW$. This finishes the proof of the claim. 
    
    Consider an orthonormal basis $\{e_1,\ldots, e_d\}$ for $T_xM$ and let $T:T_xM\to T_xM$ be the linear map with $T(e_i)=v_i$. Clearly, $|\det T|=1$,  $\hs{T}=\hs{{\bf w}}\leq H$, and for any $n\geq 0$, $\hs{\pd_xf^n\circ T}=\hs{\ef f^n({\bf w})}\leq H$. The last inequality holds since $\ef f^n  ({\bf w})\in \cW$. Using (\ref{eq:norm-1}) and (\ref{eq:norm-3}), one obtains 
    \begin{equation}\label{eq:bound-kappa}
        \|\pd_xf^n\|\leq \|\pd_xf^n \circ T\|.\|T^{-1}\|\leq \|\pd_xf^n \circ T\|.\|T\|^{d-1}\leq H^d.
    \end{equation}
    Again, by (\ref{eq:norm-3}), this implies that the sequence $\{f_i\}_{i=1}^\infty$ is $\kappa$-conformal at $x$ for $\kappa=H^{d^2}$.  
    \end{proof}

    \begin{remark}\label{rm:explicit-dep-kappa-W}
    As we observed in the proof of Theorem \ref{thm:conformality}, one can give an explicit form of dependence of $\kappa$ to $\sup_{{\bf w}\in \cW}\hs{{\bf w}}$ and vice versa. Indeed, if  item (a) in this theorem holds for some $\kappa>1$, then according to \eqref{eq:bound-norm-W}, item (b) holds for a set $\cW$ with
    \(\sup_{{\bf w}\in \cW}\hs{{\bf w}}\leq \kappa\sqrt{d}\). \\Meanwhile, if item (b) holds for some non-empty fiberwise-bounded set $\cW\subset \cE(M) $, then by \eqref{eq:bound-kappa} and \eqref{eq:norm-3}, item (a) will be true for  
    \[\kappa=(\sup_{{\bf w}\in \cW}\hs{{\bf w}})^{d^2}.\]
    \end{remark}

    \begin{corollary}[Stable quasi-conformality]\label{cor:st-QC}
    Let $V\subseteq M$ be an open set and $\cF\subseteq \diff^1_{\rm loc}(M)$. Assume that $\cW$ is an open subset of $\cE(M)$ with compact closure satisfying $\pi(\cW)=V$ and
    \begin{equation}\label{eq:cover-W}
    \overline{\cW}\subseteq \bigcup\limits_{f\in \cF}(\ef f)^{-1}(\cW).
    \end{equation}
    Then, there exists $\kappa>1$, depending only on $\cW$, such that for any family $\tilde{\cF}$
     sufficiently close to $\cF$ in the $\cC^1$ topology,
    $\IFS(\tilde{\cF}\loc{V})$ has a $\kappa$-conformal orbit-branch at every point of $V$.
    \end{corollary}
    \begin{proof}
    By the compactness of $\overline{\cW}$ and the openness of $\cW$, the same covering property (\ref{eq:cover-W}) holds for any family $\tilde{\cF}$
     sufficiently close to $\cF$ in the $\cC^1$ topology.
   Hence, the conclusion follows from Theorem \ref{thm:conformality} and Remark \ref{rm:explicit-dep-kappa-W}  for $\kappa=(\sup_{{\bf w}\in \cW}\hs{{\bf w}})^{d^2}$ which depends only  on $\cW$. 
    \end{proof}
    
    \begin{figure}[t] 
        \centering
         \includegraphics[width=.4\textwidth]{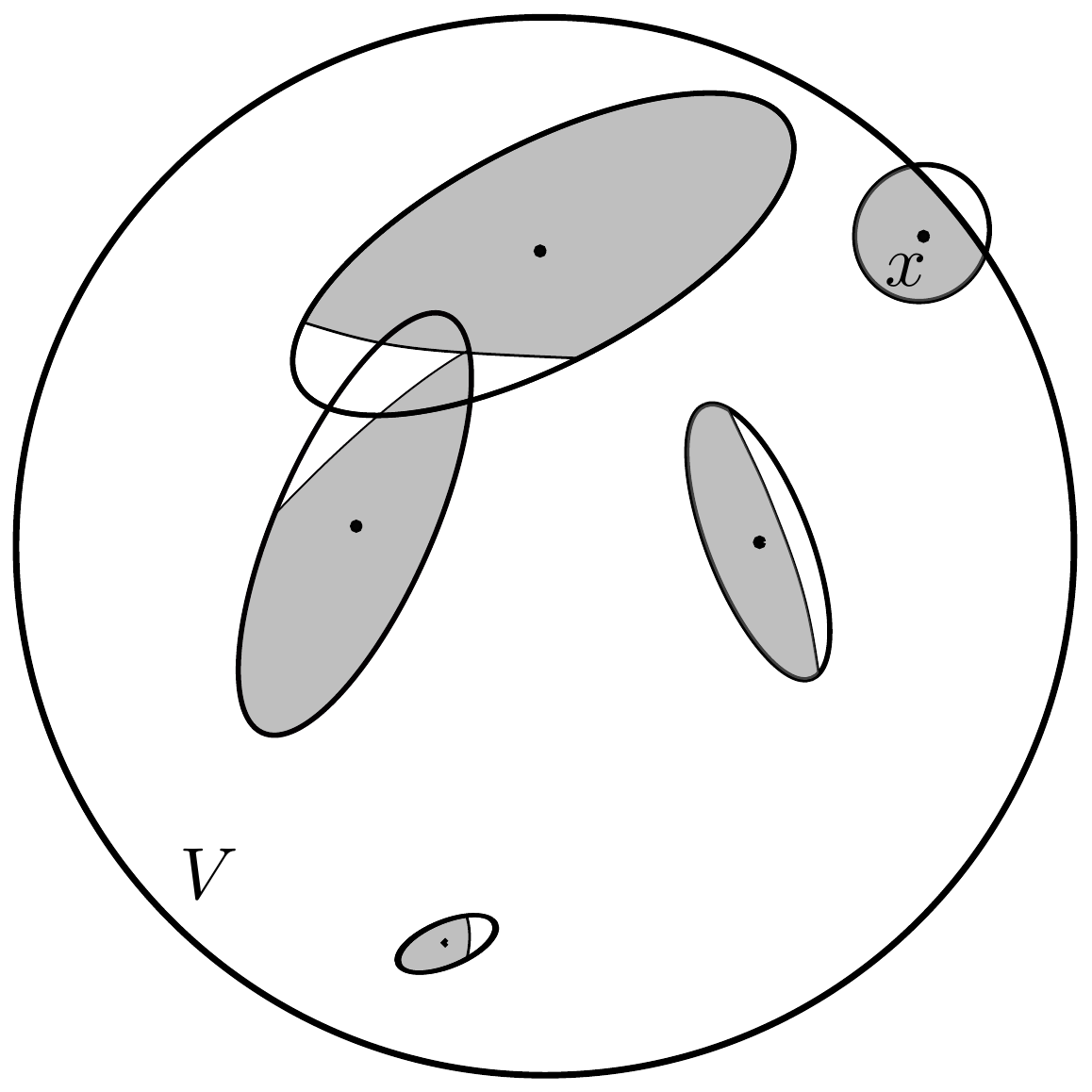}
         \caption{Covering condition (\ref{eq:cover-W}): diversity of maps corresponding to  $x\in V$, i.e. maps in the subfamily $\cF_x$.}
         \label{fig:stable-qc}
    \end{figure}

    \subsubsection{Reinterpretation of the covering property}

    The covering condition (\ref{eq:cover-W}) says that for every $x\in \overline{V}$ and ${\bf w}\in \cW_x:=\pi^{-1}(x)\cap \cW$,  there exists $f\in \cF$ such that $f(x)\in V$ and $\ef f ({\bf w})\in \pi^{-1}(f(x))\cap \cW$. 

    Another way to see this condition is the following. For a point $x$, denote by $\cF_x$ the set of all elements $f\in \cF$ with $x\in \mathrm{Dom}(f\loc{V})$. Then, the covering condition \eqref{eq:cover-W}  is equivalent to 
    \begin{equation}\label{eq:covering-fiber}
        \overline{\cW_x}\subseteq \bigcup_{f\in \cF_x} (\ef f)^{-1}(\cW_{f(x)}),
    \end{equation}
    for every $x\in \overline{V}$. 
    Roughly speaking, this means that over every point $x$ there are several maps in the family $\cF$ with diverse directions of contraction and expansion for the normalized derivative that allows obtaining the covering  condition \eqref{eq:covering-fiber}. This yields the quasi-conformality along an orbit-branch of $x$ (see Figure \ref{fig:stable-qc}).

    \subsection{Covering condition for $\sldr$ actions}

     In this subsection, we aim to find a local method to provide families of maps satisfying the covering condition (\ref{eq:cover-W}). Smooth action localized to small open balls in manifolds can be transformed to actions on open subsets of $\R^d$. Therefore, we may explore (\ref{eq:cover-W}) when, $M=\R^d$ (or an open subset of $\R^d$). In this case,  $\cE(\R^d)$ is isomorphic to the trivial fiber bundle $\R^d\times \pmsldr$ and the action of  $\ef f$ on fibers in nothing but the product of matrices in $\pmsldr$.  More precisely, for $f\in\diff^1_{\rm loc}(\R^d)$, $\ef f$ maps $(x,A)\in \R^d\times \pmsldr$ to $(f(x),  A_x A)$ where $A_x:=\pd_x f\in \pmsldr$.  In other words, the covering condition (\ref{eq:cover-W}) will be reduced to finding a map $f\in \cF$ such that $A_xA$  is in the bounded set $\cW_{f(x)}$.
    Observe that $\R^d\times \sldr$ is invariant under $\pd_x f$  for an orientation-preserving $f\in \diffloc{1}(\R^d)$. In particular, it is enough to satisfy the covering condition for $\cW\subseteq \R^d\times \sldr$. So,  we investigate the covering property for the action of $\sldr$ on itself.
    This subsection provide a method of doing that and its main result (Lemma \ref{lem:sldr-covering}) will be used in the proof of Theorem \ref{thm:example}.
\medskip

    Recall that the sequence  $\{D_i\}_{i=1}^\infty$ in $\sldr$ is quasi-conformal if and only if  the set $\{D_n\cdots D_1{\gap}n\in \N\}$  has compact closure. Also, the sequence $\{D_i\}_{i=1}^\infty$ is $\kappa$-conformal, if for any $n\in \N$, $D_n\cdots D_1$ is $\kappa$-conformal.\\ 
    
    \begin{figure}[h]
       \centering
       \includegraphics[width=.9\textwidth]{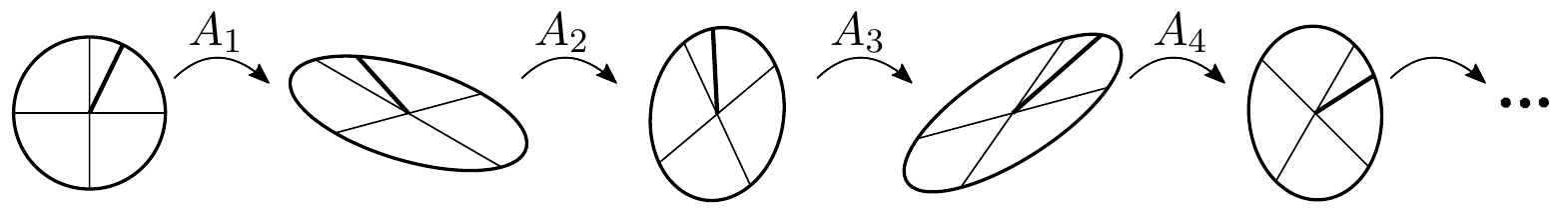}
       
       \caption{A quasi-conformal sequence of matrices}
    \end{figure} 
    
    \begin{definition}
    For $\kappa\geq 1$ and $\cD\subseteq \sldr$, we say $\genplus{D}$ has a $\kappa$-conformal branch, if there exists a $\kappa$-conformal sequence in $\cD$.
    In addition, when $\cD$ is finite, we say it robustly has $\kappa$-conformal branches, if for every $\tilde{\cD}$
    sufficiently close to 
    $\cD$, $\langle{\tilde{\cD}}\rangle^+$ has $\kappa$-conformal branches.
    \end{definition}

    \subsubsection{Behaviours of typical branches.}
    
    For a finite subset $\cD\subseteq \sldr$, if  the closure of $\genplus{D}$ is not compact, typical branches are not quasi-conformal in a probabilistic sense. More precisely, by assigning positive probabilities to the elements of $\cD$, almost every branch with respect to the product measure on $\cD^\N$ is not $\kappa$-conformal for any $\kappa>1$. 
Indeed, almost every branch in $\cD^\N$  contains all the elements of $\genplus{D}$ as sub-words, and all the sub-words of a $\kappa$-conformal branch are $\kappa^2$-conformal (see \S \ref{subsec:qc-def}). 
This implies that if the set of $\kappa$-conformal branches has positive measure (for some $\kappa$), then all the elements of $\genplus{D}$ are $\kappa^2$-conformal. This is impossible if the closure of $\genplus{D}$ is not compact due to \eqref{eq:norm-4}.

    Also, this is worth mentioning that for fixed $n\in \N$, for almost every $n$-tuple $\cD\in \sldr^n$ (w.r.t. the natural measure), the Lyapunov spectrum associated with the random product of elements of $\cD$ is non-degenerate provided that we assign positive weights to the elements of $\cD$. In such cases,
    for almost every branch with respect to the product measure on $\cD^\N$, the norm of the products diverges exponentially to infinity (cf. \cite{Viana_2014}).

    \subsubsection{Existence of bounded branches}
    The following lemma, which is an analogue of Theorem \ref{thm:conformality} for the action of $\sldr$  on itself, expresses that the covering condition leads to the existence of bounded branches for a finitely generated semigroup.      Moreover, Corollary \ref{cor:bounded-branch} for this special case evidently implies that given $\kappa>1$, any open neighbourhood $\cU$ of the identity in $\sldr$ with compact closure has a finite subset $\cD$ with robustly $\kappa$-conformal branches.
    
    \begin{lemma}\label{lem:U-branch}
    Let $\cD\subseteq \sldr$ be a finite set. Then, there exists $\kappa>1$ such that $\genplus{D}$ has a  $\kappa$-conformal branch if and only if $ ~\overline{\cU}\subseteq \cD^{-1}\cU$ for some $\cU\subseteq \sldr$ with compact closure. Moreover, if $\cU$ is open,  $\genplus{\cD}$ robustly has a $\kappa$-conformal branch for some $\kappa>1$.    
    \end{lemma}
    \begin{proof}
      Let $V=\R^d$ and $\cW=\R^d \times \cU$. Considering the natural action of $\sldr$ on $\R^d$, $\cD$ can be seen as a family in $\diff(\R^d)$. Then, the first part of the lemma follows from Theorem \ref{thm:conformality}. The second part is similar. Note that $\overline{\cU}\subseteq \cD^{-1}\cU$ implies that for any family $\tilde{\cD}$ sufficiently close to $\cD$, ${\cU}\subseteq \tilde{\cD}^{-1}{\cU}$ holds and so the second conclusion is again a consequence of Theorem \ref{thm:conformality}.
    \end{proof}
    
    We can deduce the following corollary from Lemma \ref{lem:U-branch}.
    
    \begin{corollary}\label{cor:bounded-branch}
    For any $\kappa>1$ and any open set $\cU\subseteq \sldr$ containing the identity, there is a finite set $\cD\subseteq \cU$ such that $\genplus{D}$ robustly has a $\kappa$-conformal branch. 
    \end{corollary}
    \begin{proof}
    Consider small open neighbourhood $\cV$ of the identity with compact closure such that every element of $\overline{\cV}$ is $\kappa$-conformal. Clearly, $\overline{\cV}\subseteq \cU^{-1}\cV=\bigcup_{u\in \cU}u^{-1}\cV$.  By the compactness of $\overline{\cV}$, one can choose a finite set $\cD \subseteq \cU$ with $\overline{\cV}\subseteq \cD^{-1}\cV$. Thus, the conclusion follows from Lemma \ref{lem:U-branch}.
    \end{proof}
      
    \begin{remark}
    Lemma \ref{lem:U-branch} and Corollary  \ref{cor:bounded-branch} can be stated for the existence of bounded branches in an abstract setting for more general topological groups. However, due to the applications for the derivatives of smooth maps, in this paper, the discussion is restricted to the special cases of $\sldr$ and $\R^d$.      
    \end{remark}

    \subsection{Sufficient conditions for covering:  explicit construction}
    Corollary \ref{cor:bounded-branch} is existential and does not introduce elements of  $\cD$ explicitly and does not even give any estimate for the cardinality of this set. The next lemma guarantees the covering of small open sets with $d^2$ elements.

    \begin{lemma}\label{lem:sldr-covering}
    For any neighbourhood $\cU_0$ of the identity in $\sldr$, there exists an open set $\cU\subseteq \cU_0$ and a finite set $\cD\subseteq \cU_0$ with  $d^2$ elements such that $\overline{\cU}\subseteq \cD^{-1}\cU$.
    \end{lemma}
    
    Denote by $\mathfrak{sl}(d,\R)$ the Lie algebra of $\sldr$ which consists of all $d\times d$ real matrices whose traces are equal to zero. Furthermore, here $\matrixexp$ denotes the exponential function from a neighbourhood of the zero matrix in $\mathfrak{sl}(d,\R)$ to a neighbourhood of the identity in $\sldr$  which verifies the Baker-Campbell-Hausdorff formula.  
    It will be used in the proof of the next lemma.
    \begin{lemma}\label{lem:lie-algebra}
    For any neighbourhood $\mathfrak{U}_0$ of the zero matrix in $\mathfrak{sl}(d,\R)$, there exist an open subset $\mathfrak{U}$ and a finite  subset $\mathfrak{D}=\{w_1,\ldots, w_{d^2}\}$ of  $\mathfrak{U}_0$ such that $\overline{\matrixexp(\mathfrak{U})} \subseteq \matrixexp(\mathfrak{D})^{-1} \matrixexp(\mathfrak{U})$.
    \end{lemma}
    For the proof, we use the following notation. For a connected open set $U\subseteq \R^N$ and  small $t>0$, denote
    \begin{equation}\label{Notation:U_alpha}
    U_t:=\{x\in U{\gap} d(x,\partial U)>t\}.
    \end{equation}
    In addition, the following observation will be used in the proof of Lemma \ref{lem:lie-algebra}. 
    \begin{lemma}\label{lem:image_of_close_diffeos}
    Let $U\subseteq \R^N$ be an open set with $\overline{U}$ homeomorphic to the closed unit disk. Suppose that   $\varphi_0,\varphi_1$ are two continuous maps defined on a neighbourhood of $\overline{U}$ and are homeomorphisms onto their images. If there exists $t>0$ such that for any $x\in \overline{U}$, $|\varphi_0(x)-\varphi_1(x)|<t$, then $\big( \varphi_0(U)\big)_t\subseteq \varphi_1(U)$. 
    \end{lemma}
    
    The proof of this lemma is based on considering an affine homotopy between $\varphi_0\big|_{\partial U}$ and $\varphi_1\big|_{\partial U}$, then showing that the image of homotopy does not intersect $\big(\phi_0(U)\big)_t$. Further details are left to the reader.
    \begin{proof}[Proof of Lemma \ref{lem:lie-algebra}]
    Suppose $v_1,\ldots, v_{N}$ are $N$ points in $\R^{N-1}$ with $|v_i|=1$ such that the origin is contained in their convex hull. Denote the interior of their convex hull by $\Delta$. Clearly, for any small positive $t\in \R$, $\overline{\Delta}\subseteq \bigcup_{i=1}^N (\Delta+t v_i)$. Since $\overline{\Delta}$ is compact and $\Delta+t v_j$'s are all open, one can find sufficiently small $c>0$ such that $\overline{\Delta}\subseteq \bigcup_{i=1}^N(\Delta+t v_i)_{c}$ (following the notation introduced in (\ref{Notation:U_alpha})). As the whole construction is invariant under homothety, for any  $r>0$,
    \begin{equation}\label{eq:Linear_Convering_of_Simlex_Delta}
    r\overline{\Delta}=\overline{r\Delta}\subseteq \bigcup_{i=1}^N (r\Delta+tr v_i)_{cr}.
    \end{equation}
    
    By identification of  $\mathfrak{sl}(d,\R)$ with $\R^{N-1}$ for $N=d^2$,  $v_i$'s can be seen as $d\times d $ matrices with zero trace. 
    For any $1\leq j\leq N$, and any sufficiently small  $r>0$, we define $u^{(r)}_j:r\Delta\to \mathfrak{sl}(d,\R)$ as the following 
    $$u_j^{(r)}(x)=\matrixexp^{-1}\big(\matrixexp(trv_j)\matrixexp(x)\big).$$
    By the Baker-Campbell-Hausdorff formula for the Lie groups (see for instance \cite{Rossmann}),
    \[\matrixexp^{-1}\big(\matrixexp(trv_j)\matrixexp(x)\big)=x+trv_j+\varepsilon_j(x,tr),\]
    where $\varepsilon_j$ satisfies $|\varepsilon_j(x,s)|<E_j|x|s$ for some  $E_j>0$.
    Thus, whenever $r<\min\limits_{1\leq j \leq N}\{\frac{c}{E_jt}\}$, for $1\leq j\leq N$ one has,
    \[\sup_{x\in\overline{r\Delta}}\big|u_j^{(r)}(x)-(x+trv_j)\big|=\sup_{x\in\overline{r\Delta}}|\varepsilon_j(x,rt)|\leq E_jrt\sup\limits_{x\in\overline{r\Delta}}|x|\leq E_jr^2t<cr.\]
    By  Lemma \ref{lem:image_of_close_diffeos}, for every $1\leq j \leq N$, $r\Delta+trv_j\subseteq u_j^{(r)}(r\Delta)$ and so by     (\ref{eq:Linear_Convering_of_Simlex_Delta}),
    $$ \overline{r\Delta}\subseteq \bigcup\limits_{j=1}^N u^{(r)}_j(r\Delta).$$
    Finally, as the exponential map is a diffeomorphism on $u^{(r)}_j(r\Delta)$ and on $r\Delta$, 
    \begin{equation}\label{eq:covering_equation}
    \overline{\matrixexp(r\Delta)}=\matrixexp(\overline{r\Delta})\subseteq \matrixexp\Big(\bigcup\limits_{j=1}^N u_j^{(r)}(r\Delta)\Big) =\bigcup\limits_{j=1}^N\Big(\matrixexp\big(u_j^{(r)}(r\Delta)\big)\Big).
    \end{equation}
    Meanwhile, by the definition of $u_j^{(r)}$, we have $\matrixexp(u_j^{(r)}(r\Delta))=\matrixexp(trv_j)\matrixexp(r\Delta).$
    Thus, the conclusion of Lemma \ref{lem:lie-algebra}  follows from (\ref{eq:covering_equation}) by taking
    $w_j:=-trv_j$ and $\mathfrak{U}:=r\Delta$ {for sufficiently small $r$.}
    \end{proof}
    
    \begin{remark}\label{rm:minimal-blender}
    {One can start with any simplex $\Delta$ in $\mathfrak{sl}(d,\R)\simeq \R^{d^2-1}$ containing the origin in the interior to provide an explicit formula for $\mathfrak{D}$ in Lemma \ref{lem:lie-algebra}}.
    \end{remark}
    
    \begin{proof}[Proof of Lemma \ref{lem:sldr-covering}]
    Consider an open neighbourhood $\mathfrak{U}_0$ of the zero matrix in $\mathfrak{sl}(d,\R)$ such that the $\matrixexp$ function is a diffeomorphism in a neighbourhood of  $\overline{\mathfrak{U}_0}$ and $\exp(\overline{\mathfrak{U}_0})\subseteq \cU_0$. 
    Now, applying Lemma \ref{lem:lie-algebra} to $\mathfrak{U}_0$, one can get $\mathfrak{U},\mathfrak{D}$. Then, $\cU:=\matrixexp(\mathfrak{U})$ and $\cD:=\matrixexp(\mathfrak{D})$ satisfy the conditions of Lemma \ref{lem:sldr-covering} and the proof is finished.
    \end{proof}
   
    \subsection{Sufficient conditions for covering:  geometric construction} 
   
    This subsection states another approach to derive a sufficient condition leading to the covering property with respect to a fiberwise-bounded subset of $\cE(M)$.  It will be used in the proof of some of the main results, including Theorems \ref{thm:any-dim-sullivan} and \ref{thm:sphere-2}.
    In comparison to Lemma \ref{lem:sldr-covering}, the next lemma is more geometric and does not assume the generators being close to the identity. See also Question \ref{q:covering}.

    \begin{lemma}\label{lem:direc-cont}
    Let $\cF\subseteq \diff_{\rm loc}^{1}(M)$ and $V\subseteq M$ be an open set with compact closure. Assume that for any $(x,v)\in T^1M$ with $x\in \overline{V}$, there exists $f\in  \cF$ satisfying $f(x)\in V$ and $\|\pd_x f|_{v^\perp}\|<1$. Then, there exists an open set $\cW \subseteq \cE(M)$ with compact closure such that $\pi(\cW)=V$ and $\overline{\cW} \subseteq  \bigcup_{f\in \cF} (\ef f)^{-1}(\cW)$.
    \end{lemma}

    \begin{proof}
    Since the set $\{(x,v)\in T^1M{\gap} x\in \overline{V}\}$ is compact, one can find a finite subset $\cF_0\subseteq \cF$ and $\epsilon>0$ such that for any $(x,v)\in T^1M$ with $x\in \overline{V}$, there exists $f\in \cF_0$ satisfying $f(x)\in V$ and $\|\pd_x f|_{v^\perp}\|<1-\epsilon$.
    
    Denote $\Theta:=\max\{\|\pd_xf\|{\gap}f\in \cF_0, x\in \overline{V}\}$. Let $H\geq \Theta^4/\epsilon$
    be a real number. Then, consider $\cW\subseteq \cE(M)$, defined by  
    \begin{equation}\label{eq:def-W}
        \cW:=\{{\bf w}\in \cE(M){\gap}\pi({\bf w})\in V~ \text{and} ~\hs{{\bf w}}< H\}.
    \end{equation}
    Clearly, $\cW$ is an open subset with compact closure and $\pi(\cW)=V$. It is enough to prove that (\ref{eq:cover-W}) holds for the subfamily $\cF_0$ and  $\cW$ defined by (\ref{eq:def-W}). Consider ${\bf w}=(x,(v_1,\ldots,v_d))\in \overline{\cW}$. So, $\hs{{\bf w}}\leq H$. If $\hs{{\bf w}}<\Theta^{-1}H$ and $f\in \cF_0$ with $f(x)\in V$, then by (\ref{eq:norm-2}),
    \[\hs{\ef f ({\bf w})}=\Big(\sum\limits_{i=1}^d |\pd_x f (v_i)|^2\Big)^\frac{1}{2}\leq \|\pd_x f\|.\hs{{\bf w}}< H.\]
    Thus, $\ef f({\bf w})\in \cW$ and consequently, ${\bf w}\in \bigcup_{f\in \cF_0} (\ef f)^{-1}(\cW)$. Now,  assume that  $\hs{{\bf w}}\in [\Theta^{-1}H, H]$. Without loss of generality, assume that $|v_1|\geq \cdots \geq |v_d|$. Clearly, $|v_d|\leq \sqrt{H/d}$. So, 
    \[|v_d|\Theta \sqrt{d/\epsilon}\leq \Theta \sqrt{H/\epsilon} \leq \Theta^{-1}H\leq \hs{{\bf w}}=\sqrt{|v_1|^2+\cdots +|v_d|^2}\leq \sqrt{d}|v_1|.\]
    This implies that $\Theta^2 |v_d|^2\leq \epsilon |v_1|^2$. Let $v$ be a vector perpendicular to $v_1,\ldots,v_{d-1}$. By assumption, one can choose $f\in \cF_0$ with $f(x)\in V$ and $\|\pd_x f|_{v^\perp}\|<1-\epsilon$. Then, 
    \begin{align*}
        \hs{\ef f ({\bf w})}^2&=\sum\limits_{i=1}^d |\pd_x f (v_i)|^2\\
        &<(1-\epsilon)^2 (|v_1|^2+\cdots +|v_{d-1}|^2)+\Theta^2|v_d|^2\\
        &\leq (1-\epsilon)^2(|v_1|^2+\cdots+|v_{d-1}|^2)+\epsilon |v_1|^2<\hs{{\bf w}}^2\leq H^2.
        \label{eq:bound-HS}
    \end{align*}
    Therefore, ${\bf w}\in \bigcup_{f\in \cF_0} (\ef f)^{-1}(\cW)$. 
    \end{proof}
    
    \section{Expanding sequences} \label{sec:expanding}
    Here, we prove two technical results which will be used for showing ergodicity of quasi-conformal blenders. The first one provides a precise control of geometry under quasi-conformal expanding sequences. The second one is the standard bounded distortion lemma adapted to our setting of expanding sequences of local maps.

    Throughout the section,  $M$ is a closed manifold of dimension $d$. Recall that for a sequence $\{f_i\}_{i=1}^\infty$ in  $\diff_{\rm loc}^{1+\alpha}(M)$, for $i\in \N$ $f^i:=f_i\circ \cdots \circ f_1$, $f^{-i}:=(f^i)^{-1}$, and $f^0:=\mathrm{Id}$.  We denote the open ball of the radius of $r>0$ around the origin in $\R^d$ with $\Eucball{r}$.

    \subsection{Control of geometry}
     Our main goal in this subsection is to prove the following theorem. 
     
    \begin{theorem} \label{thm:control-shape-M}
    Let  $\{f_i\}_{i=1}^\infty$ be a sequence in $\diff_{\rm loc}^{1+\alpha}(M)$
    with bounded $\cC^{1+\alpha}$ norm. 
    Let also $x\in M$ and $\rho>0$ be such that for any $n\in \N$, $f_n$ is defined on $\ball{f^{n-1}(x)}{\rho}$. If the sequence is quasi-conformal expanding at $x$, then there exist {$\xi_0>0$} and $\theta>1$ such that for any {$\xi\in (0,\xi_0]$} and $n\in \mathbb{N}$,    
     \[ \ball{x}{r_n}\subseteq  f^{-n}(\ball{f^n(x)}{{\xi}})\subseteq     \ball{x}{\theta r_n},\]
    for some $r_n>0$.
    \end{theorem}
    
    In other words, one obtains a control of geometry of the iterations of a ball from certain assumptions on the derivatives at the center. The passage from linear to nonlinear follows from precise estimates on pseudo-orbits of corresponding product of matrices. The proof of this theorem occupies the entire subsection. As it is a local statement, we prove Theorem \ref{thm:control-shape-M} by showing similar statements on Euclidean space and for uniformly contracting sequences of local diffeomorphisms. We also prove a stronger version of Theorem \ref{thm:control-shape-M} in Theorem \ref{thm:st-control-shape}.

\medskip

    Given $R,C>0,\kappa\geq  1>\uppercont>\lowercont>0 $ and $\alpha\in (0,1)$. For $N\in \N\cup \{\infty\}$, we consider  the following hypotheses for the sequence $\{h_n\}_{n=1}^N$.
    
    \medskip
    \begin{itemize}
        \item[(H0)]   $h_n:\Eucball{R}\to h_n(\Eucball{R})$ is a $\cC^{1+\alpha}$ diffeomorphism fixing the origin,
        \item[(H1)] $\big\|h_n\big\|_{\cC^{1+\alpha}}< C$, 
        \item[(H2)] for any $y\in \Eucball{R}$, 
        $\lowercont<m(D_y h_n)\leq \|D_yh_n\|<\uppercont$,
        \item[(H3)] $h^n$ is $\kappa$-conformal at the origin.
    \end{itemize}
    \medskip
    Note that, it follows from (H2) that $h_n(\overline{\Eucball{R}})\subseteq \Eucball{R}$. 
    \begin{figure}[h]
    \includegraphics[width=.95\textwidth]{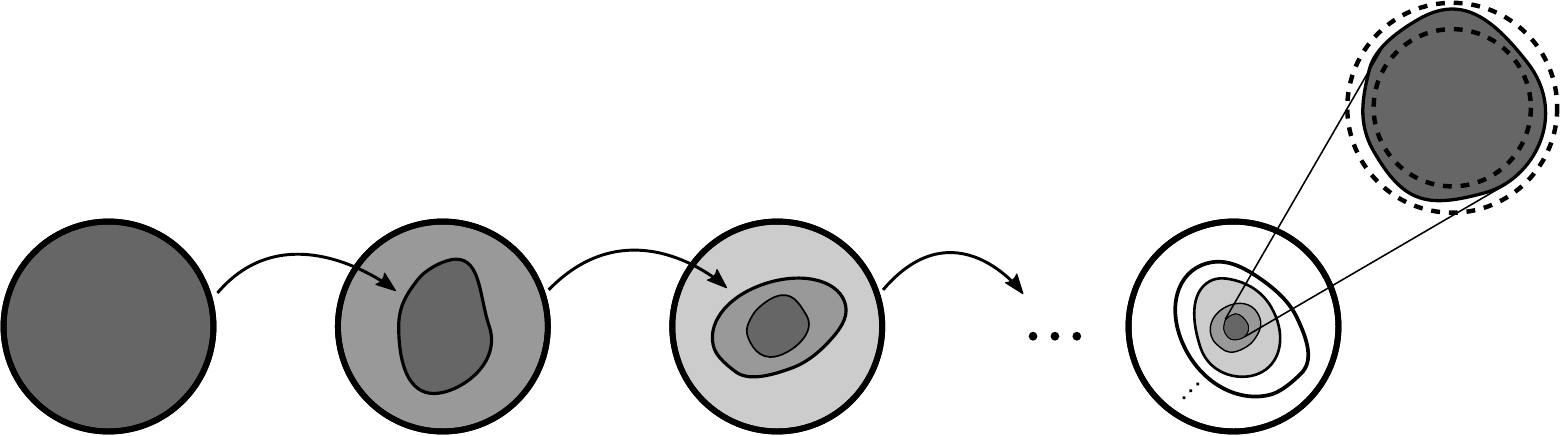}
    \caption{Sequence of maps satisfying (H0)-(H3) in Lemma \ref{lem:iteration-balls-Rd}, where the iterations of small ball remain almost round.}
     \label{fig:iteration-balls}
     \end{figure}
    \begin{lemma}\label{lem:iteration-balls-Rd}
     Let $R,C>0$, $\kappa\geq 1 >\uppercont> \lowercont>0$ and  $\alpha\in (0,1)$. Then, there exist  $\xi_0\in (0,R)$ and ${\gamma},\theta>1$ such that for any $\xi\in {(0,\xi_0]}$, any $n\in \N$ and any sequence $\{h_j\}_{j=1}^n$ of maps satisfying (H0)-(H3),
    \begin{equation}\label{eq:weak-lin-to-nonlin}
        \Eucball{r/\theta}\subseteq h^n(\Eucball{\xi})\subseteq \Eucball{\theta r},
    \end{equation}
    for $r=\xi|\det D_0h^n|^{1/d}$. Moreover, for any $|x|\leq \xi$,
    \begin{equation}\label{eq:st-lin-to-nonlin}
        |h^n(x)- D_0h^n(x)|<{\gamma}|\det D_0h^n |^{1/d}|x|^{1+\alpha}.    
     \end{equation}
    \end{lemma}

    This lemma follows from the next one on sequences of linear maps satisfying conditions (H2)-(H3). It establishes precise estimates on the difference between the orbits and certain pseudo-orbits.  Indeed, the following lemma plays a crucial role in the proof of all the main results.

    \begin{lemma}[Key Lemma]\label{lem:pseudo-orbit}
    Let $C>0$,  $\kappa\geq 1>\uppercont>\lowercont>0$ and $\alpha\in(0,1)$. Then, there exist $\ettainemma,\gamma>0$ such that for any $n\in \N$, any sequence  $\{D_i\}_{i=1}^n$ of matrices in $\mathrm{GL}(d,\R)$ and any sequence $\{y_i\}_{i=0}^n$ in $\R^d$ satisfying 
    \begin{itemize}
    \item[(C1)] Contraction:  $\lowercont\leq m(D_i)\leq \|D_i\|\leq \uppercont<1$,
    \item[(C2)] Quasi-conformality: $\|D_{i}D_{i-1}\cdots D_1 \|\leq \kappa \cdot m(D_{i}D_{i-1}\cdots D_1)$,
    \item[(C3)] $|y_0|\leq \ettainemma$ and $|y_{i}-D_{i}y_{i-1}|<C|y_{i-1}|^{1+\alpha},$
    \end{itemize}
  for any $1\leq i \leq n$, the following holds \begin{equation}\label{eq:key-lem-result}
        |y_n-D_n\cdots D_1y_0|<\gamma  |\det (D_n\cdots D_1)|^{1/d} |y_0|^{1+\alpha}.
    \end{equation}
    \end{lemma} 
    \begin{proof}
    To simplify the notation, we denote $D_{j,i}:=D_{j}D_{j-1}\cdots D_i$ for $ j>i$.
    \medskip 
    
    The proof has several steps. 
    First we prove a useful estimate for the determinant of large blocks of products (of some length $K$) in terms of their norm (Claim 1). So, it would be convenient to replace the sequence $\{D_i\}_i$ with its blocks of products $\{D_{(i+1)K, iK+1}\}_i$. Then Claim 2 assures that this new sequence also satisfies the assumptions of the lemma, particularly 
    (C3). The next two steps prove the lemma for this sequence. 
    In other words, Claims 1-4 prove the conclusion of Lemma \ref{lem:pseudo-orbit} when $n$ is divisible by $K$. Finally, Claim 5 proves the Key Lemma by extending it to all $n$.
    
        \medskip

   Since $\{D_i\}$ satisfies  (C2), we conclude that for any $j>i$, $D_{j,i}$ is $\kappa^2$-conformal. By \eqref{eq:ineq-norm-conorm}
    \[\frac{\|D_{j,i}\|}{|\det D_{j,i}|^{1/d}}\leq \frac{\|D_{j,i}\|}{m(D_{j,i})}\leq \kappa^2,\text{and}~ \frac{|\det D_{j,i}|^{1/d}}{m(D_{j,i})}\leq \frac{\|D_{j,i}\|}{m(D_{j,i})}\leq \kappa^2, \]
   
    So the condition (C2) implies the following.
   
    \begin{itemize}
       \item[(C2)$^\prime$] For any $j> i \geq 1$,  $\uzetatwo |\det D_{j,i}|^{1/d}\leq m(D_{j,i})\leq \|D_{j,i}\|\leq \ozetatwo |\det D_{j,i}|^{1/d}$.
    \end{itemize}
    
    \medskip

    \noindent
    \textbf{Claim 1.}\emph{
    For any $\alpha'>0$, there exist $K=K(\kappa,d,\uppercont,\lowercont,\alpha')\in \N$ and $\tau>0$ such that for any $i\geq 0$,
    \begin{equation}\label{eq:good-block}
    \big(\|D_{i+K,i+1}\|+\tau\big)^{1+\alpha'}\leq |\det D_{i+K,i+1}|^{1/d}.  
    \end{equation}
    }
    \begin{proof}
    For $t\in \R^+$, denote $\varphi(t):=t^\frac{1}{1+\alpha'}-\ozetatwo t$. Since $1+\alpha'>1$, there exists $T=T(\alpha',\kappa)>0$, such that  $\varphi(t)$ is positive and increasing on $(0,T)$. Let $K\in \N$ be large enough such that $\uppercont^K<T$ and $\tau:=\varphi(\uzetatwo\lowercont^K)$. Now, the conclusion easily follows from (C2)$^\prime$ and (C1).
    \end{proof}
    
    \medskip
    
    Note that from 
    (C3), $|y_{i+1}|<C|y_i|^{1+\alpha}+\|D_{i+1}\||y_i|$.  Therefore, if $|y_0|\leq \xi':=C^{-\frac{1}{\alpha}}(1-\uppercont)^\frac{1}{\alpha}$, then for any $i\geq 0$, $C|y_i|^\alpha+\|D_{i+1}\|\leq 1$ and so
    \begin{equation}\label{eq:y_i-dec}
        |y_{i+1}|\leq |y_i|.
    \end{equation}
     Define $\epsilon_i:=y_{i+1}-D_{i+1}y_{i}\in\R^d$. For any pair $i,k\geq 0$, one obtains an explicit formula for $y_{i+k}$ in terms of $y_{i}$ and $\{\epsilon_j{\gap}i\leq j< i+k\}$, 
    \begin{align}
    y_{i+k}   & = D_{i+k}y_{i+k-1}+\epsilon_{i+k-1}=D_{i+k}(D_{i+k-1} y_{i+k-2} +\epsilon_{i+k-2})+\epsilon_{i+k-1}=\cdots\nonumber\\
        & =D_{i+k,i+1}y_i+ \sum\limits_{j=i+2}^{i+k}D_{i+k,j}\epsilon_{j-2}+\epsilon_{i+k-1}.\label{eq:explicit formula for y_0}
   \end{align}
   
   \medskip
   \noindent
   \textbf{Claim 2.} \emph{If $|y_0|\leq \xi'$, then for any $i\geq 0$ and $k\geq 1$,
    \[|D_{i+k,i+1}y_i-y_{i+k}|<Ck|y_i|^{1+\alpha}.\]}
    \begin{proof}
    By (\ref{eq:explicit formula for y_0}) and since $\|D_j\|<1$, 
    \[|y_{i+k}-D_{i+k,i+1}y_i|\leq  \sum\limits_{j=i}^{i+k-1}|\epsilon_{j}|.\]
    
    On the other hand, 
    (C3)  and (\ref{eq:y_i-dec}) imply that for any $j=i,\ldots,i+k-1$, $|\epsilon_j|<C|y_j|^{1+\alpha}\leq C|y_i|^{1+\alpha}$. This finishes the proof of Claim 2.
    \end{proof}
    Now, let $\alpha'\in (0,\alpha)$. Suppose that $K\in \N$, $\tau>0$ are the numbers provided by Claim 1. The aim of Claims 3 and 4 is to prove the conclusion of Lemma \ref{lem:pseudo-orbit} when $n$ is divisible by $K$, and then Claim 5 completes the proof for arbitrary $n$.
   
    \medskip
    \noindent
    \textbf{Claim 3.}\emph{
    There exist $\ettainemma,\tau'>0$ with $\tau'<\min\{1-\uppercont^K,\tau\}$ and $\ettainemma\leq \xi'$ such that if $|y_i|\leq \ettainemma$, then for any $p\geq 1$,}
    \begin{equation}\label{eq:claim-3}
        |y_{pK+i}|\leq \prod_{j=0}^{p-1}\big(\|D_{(j+1)K+i,jK+i+1}\|+\tau'\big)|y_i| \leq |y_i|.
    \end{equation}
    \begin{proof}
    It suffices to prove (\ref{eq:claim-3}) for $p=1$. The general case follows immediately from induction on $p$. Let $\tau',\ettainemma >0 $ be such that $\tau'<\min\{1-\uppercont^K,\tau\}$ and $\ettainemma\leq \min\{(C^{-1}K^{-1}\tau')^\frac{1}{\alpha},\xi'\}$. Then, by Claim 2,
    \[\frac{|y_{i+K}|}{|y_i|}<\|D_{i+K,i+1}\|+CK|y_i|^\alpha\leq \|D_{i+K,i+1}\|+\tau'\leq \uppercont^K+\tau'<1.\]
    This finishes the proof of Claim 3. 
    \end{proof}
        
    \medskip
    \noindent
    \textbf{Claim 4.}\emph{
    There exists $\gamma'>0$ such that for any $p\in \N$ and $q\geq 0$,  
    \[|y_{pK+q}-D_{pK+q,q+1}y_q|<\gamma' |\det D_{pK+q,q+1}|^{1/d} |y_q|^{1+\alpha},\] 
    provided that $|y_q|\leq \xi$.}
    \begin{proof}
    Define $\alpha'':=\alpha-\alpha'>0$. For simplicity, for a fixed $q$ and $1\leq i \leq  p$, we write $D'_{i}:=D_{iK+q,(i-1)K+q+1}$, $\uppercont'_i:=\|D'_i\|$, $\lowercont'_i:=m(D'_i)$. Moreover, for $0\leq i \leq p$, let $y'_i:=y_{iK+q}$. Also, denote $\epsilon'_i:=y'_{i+1}-D'_{i+1}y'_i$. Similar to (\ref{eq:explicit formula for y_0}), one gets,
    \begin{equation*}
     y'_p= D'_{p,1} y'_0 + D_{p,2}  \epsilon'_{0} +\cdots + D'_{p,p-1} \epsilon'_{p-3}+D'_p\epsilon'_{p-2}+\epsilon'_{p-1},
      \end{equation*}
    and so, 
    \begin{equation}\label{eq:upper_bound_for_y_0}
    |y'_p-D'_{p,1} y'_0| \leq  |D'_{p,2}  \epsilon'_{0}| +\cdots +|D'_p\epsilon'_{p-2}|+|\epsilon'_{p-1}|.
    \end{equation}
    Let $\uppercont':=\uppercont^K$ and $\lowercont':=\lowercont^K$. Claim 1 and $\tau'\leq \tau$ imply
    \begin{equation}\label{eq:estimate for epsilon_i part 2}
    (\uppercont'_j+\tau')^{1+\alpha}= (\uppercont'_j+\tau')^{1+\alpha'}(\uppercont_j'+\tau')^{\alpha''}\leq \delta'_j(\uppercont'+\tau')^{\alpha''}. \end{equation}
    Hence, by Claims 2 and 3, 
        \begin{align*}
    |\epsilon'_k|\leq CK|y'_{k}|^{1+\alpha}
    & \leq CK|y'_0|^{1+\alpha}\prod_{j=1}^{k}(\uppercont'_j+\tau')^{1+\alpha}\\
    &\leq CK {(\uppercont'+\tau')}^{{k\alpha''}}\delta'_k \cdots \delta'_1|y'_0|^{1+\alpha}\label{eq:estimate for epsilon_i part 1}
    \end{align*}
    By (C2)$^\prime$, for $p\geq k+2$,
    \begin{align*}
        |D'_{p,k+2}\epsilon'_{k}|&\leq \ozeta\delta'_p\cdots \delta'_{k+2} |\epsilon'_{k}|\\
        & \leq CK\ozeta\delta'_p\cdots \delta'_{k+2}{(\uppercont'+\tau')}^{{k\alpha''}}\delta'_k \cdots \delta'_1|y'_0|^{1+\alpha}\\
        & \leq \frac{CK\ozeta} {\lowercont'} (\uppercont'+\tau')^{k\alpha''}\delta'_p\cdots \delta'_{1} |y'_0|^{1+\alpha},
    \end{align*} 
    where the last inequality follows from $\delta'_{k+1}\geq \lowercont'_{k+1}$. Now, from (\ref{eq:upper_bound_for_y_0}), 
    \[|y'_p-D'_{p,1} y'_0| \leq \frac{CK\ozeta} {\lowercont'}\Big(\sum\limits_{k=0}^{p-1}(\uppercont'+\tau')^{k\alpha''}\Big)\delta'_p\cdots \delta'_{1} |y'_0|^{1+\alpha}.\]
     On the other hand,  $\uppercont'+\tau'<1$. Consequently,  $$\sum\limits_{k=0}^{p-1}(\uppercont'+\tau')^{k\alpha''}\leq C':=\sum\limits_{k=0}^{\infty}(\uppercont'+\tau')^{k\alpha''}<\infty.$$  Therefore, taking $\gamma':=\ozeta CKC'\lowercont^{-K}$ finishes the proof.
    \end{proof}
    
    \medskip
    \noindent
    \textbf{Claim 5.}\emph{
    There exists $\gamma>0$ such that if $|y_0|\leq \ettainemma$, }
    \[|y_n-D_{n,1}y_0|<\gamma |\det D_{n,1}|^{1/d}|y_0|^{1+\alpha}.\]
    \begin{proof} Let $\delta_i:=|\det D_{i}|^{1/d}$ and write $n=pK+q$ for $0\leq q <K$. Since $|y_q|\leq \ettainemma$ from Claim 4, 
    \begin{align*}
         |y_{pK+q}-D_{pK+q,q+1}y_q|\leq \gamma'\delta_{pK+q}\cdots \delta_{q+1} |y_q|^{1+\alpha}.
    \end{align*}
    Meanwhile, from (C2)$^\prime$ and Claim 2, it follows that
    \begin{align*}
        |D_{pK+q,q+1}y_q-D_{pK+q,1}y_0|&\leq \|D_{pK+q,q+1}\|.|y_q-D_{q,1}y_0|\\
        &\leq \ozetatwo \delta_{pK+q}\cdots \delta_{q+1} Cq|y_0|^{1+\alpha}.
    \end{align*}
    Finally, by Claim 4 and since $\delta_i\geq \lowercont$, one has
        \begin{align*}
        |y_{pK+q}-D_{pK+q,1}y_0|&\leq \delta_{pK+q}\cdots \delta_{q+1} |y_0|^{1+\alpha}(\gamma'+\ozeta Cq)\\
        &\leq \delta_{pK+q}\cdots \delta_{q+1} \delta_q\cdots \delta_1 |y_0|^{1+\alpha}\lowercont^{-q}(\gamma'+\ozeta Cq).
         \end{align*}
    So, the conclusion holds for $\gamma:=\lowercont^{-K}(\gamma'+\ozeta CK)$. 
    \end{proof}
    Claim 5 completes the proof of Lemma \ref{lem:pseudo-orbit}.
    \end{proof}

    Now, we are ready to complete the proof of Lemma \ref{lem:iteration-balls-Rd} and Theorem \ref{thm:control-shape-M}.
    
    \begin{proof}[Proof of Lemma \ref{lem:iteration-balls-Rd}]
    Take $x\in \Eucball{R}$. Clearly, if the sequence $\{h_j\}_{j=1}^n$ satisfies (H0)-(H3), then the conditions of Lemma \ref{lem:pseudo-orbit} are satisfied for   $D_j:=D_0h_j$, $y_j:=h^j(x)=h_i(x_{i-1})$. So, (\ref{eq:st-lin-to-nonlin}) holds for any $|x|\leq \xi_1$. In order to prove (\ref{eq:weak-lin-to-nonlin}), note that (\ref{eq:st-lin-to-nonlin}) in particular implies that 
    \begin{equation}\label{eq:triangle-ineq-1}
        |D_0h^n(x)|-\gamma |\det D_0h^n(x)|^{1/d}|x|^{1+\alpha}<|h^n(x)|,
    \end{equation}
    and 
    \begin{equation}\label{eq:triangle-ineq-2}
        |h^n(x)|<|D_0h^n(x)|+\gamma |\det D_0h^n(x)|^{1/d}|x|^{1+\alpha}.
    \end{equation}
    On the other hand, by (C2)$^\prime$, 
    \begin{equation*}
        {\kappa^{-2}} |\det D_0h^n(x)|^{1/d} |x| \leq |D_0h^n(x)|\leq {\kappa^{2}} |\det D_0h^n(x)|^{1/d} |x|.
    \end{equation*}
    This combined with  (\ref{eq:triangle-ineq-1}) and (\ref{eq:triangle-ineq-2}), implies that
    \begin{equation*}
        {\kappa^{-2}}|x|-{\gamma}|x|^{1+\alpha}<\frac{|h^n(x)|}{|\det D_0h^n(x)|^{1/d}}<{\kappa^{2}}|x|+{\gamma}|x|^{1+\alpha}.
    \end{equation*} 
    Therefore, for any $\xi$ with ${\kappa^{-2}}-{\gamma}\xi^\alpha>0$, $\Eucball{r_1(\xi)}\subseteq h^n(\Eucball{\xi})\subseteq \Eucball{r_2(\xi)}$, provided that 
     \[r_1(\xi):= ({\kappa^{-2}}-{\gamma}\xi^{\alpha})\xi|\det D_0 h^n|^{1/d}, \quad r_2(\xi):= ({\kappa^{2}}+{\gamma}\xi^{\alpha})\xi|\det D_0 h^n|^{1/d}.\]
    To prove (\ref{eq:weak-lin-to-nonlin}), take $\theta>{\kappa^{2}}$. Then, for sufficiently small $\xi>0$, ${\kappa^{2}}+{\gamma}\xi^{\alpha}< \theta$ and ${\kappa^{-2}}-{\gamma}\xi^{\alpha}>\theta^{-1}$. 
     \end{proof}

    \begin{proof}[Proof of Theorem \ref{thm:control-shape-M}]
    {Let $R_1>0$ be smaller than the radius of injectivity of the $\exp$ function on $M$ and $x_j:=f^j(x)$. Suppose the sequence is $\kappa$-conformal $\eta$-expanding at $x_0$.
    Since the sequence has bounded $\cC^{1+\alpha}$ norm, there is $\rho'\in (0,\rho)$ and $\upperexp, \lowerexp>1$ such that for any $i\in \N$ and $y\in \ball{x_{i-1}}{\rho'}$,
    \[\lowerexp< m(D_yf_i)\leq \|D_yf_i\|<\upperexp.\]
    In fact, one has $\sup_{i\in \N}\|Df_i|_{\ball{x_{i-1}}{\rho}}\|<\infty$, and if $C>\sup_{i\in \N} \|f_i\|_{\cC^{1+\alpha}}$, by \eqref{eq:diff-conorm}
    \begin{equation*}
        |m(D_{x_{i-1}}f_i)-m(D_yf_i)|\leq \|D_{x_{i-1}}f_i-D_yf_i\|  <C|x_{i-1}-y|^\alpha.
    \end{equation*}
    So, $m(D_yf_i)>\lowerexp$ for any $y\in \ball{x_{i-1}}{\rho'}$, provided that $\rho'<(C^{-1}(\eta-\lowerexp))^\frac{1}{\alpha}$}. 
    
    For $R<R_0:=(\upperexp)^{-1}\min\{R_1,\rho'\}$ and $i\in \N$, the map $\tilde{f}_i:=\exp_{x_i}^{-1}\circ f_i \circ \exp_{x_{i-1}}$ is defined on  $\Eucball{R}\subseteq T_{x_i}(M)$ and is a diffeomorphism onto its image. 
    After an isometric identification of the tangent spaces with $\R^d$, one can consider the sequence $\{\tilde{f}_i\}_{i=1}^\infty$ as a sequence of expanding maps defined on $\Eucball{R}\subseteq \R^d$. By uniform expansion of the maps, $\Eucball{R}\subseteq \tilde{f}_i(\Eucball{R})$.
    
    Next, fix $n\in \N$. The sequence $\{h_j\}_{j=1}^n$ defined by $h_j:=\tilde{f}_{n+1-j}^{-1}\mid_{\Eucball{R}}$ satisfies hypotheses (H0)-(H3) with constants independent of the choice of $n$. So, by Lemma \ref{lem:iteration-balls-Rd}, there are $\theta>1$ and $\xi_0>0$ such that for any $\xi\leq \xi_0$ and for some $r_n>0$,
    $$(h^n)^{-1}(\Eucball{r_n})\subseteq \Eucball{\xi}\subseteq (h^n)^{-1}(\Eucball{\theta r_n})\subseteq \Eucball{R}.$$
    This finishes the proof, since $(h^n)^{-1}=\exp^{-1}_{x_n}\circ f^n \circ \exp_{x_0}$ and for small $r>0$, the function $\exp_x:T_xM\to M$, maps $\Eucball{r}\subseteq T_xM$ to $\ball{x}{r}\subseteq M$.
    \end{proof}
    \subsection{Bounded distortion}
    In this subsection, we present a well-known bounded distortion lemma for a sequence of contracting maps, which permits us to control the growth of measure of iterations of measurable sets. For the purpose of completeness, we present the proof here, which is an adaptation of the classical argument to our setting. Let $R_1>0$ be smaller than the radius of injectivity of the exponential map on $M$.
  
    \begin{lemma}\label{lem:distortion-part1}
    Let   $\alpha,\lambda\in (0,1)$ and $C>0$. Then, there exists $L>1$ such that for any $R<R_1$, any $n\in \N$, any sequence $\{x_j\}_{j=0}^n$ in $M$, and any sequence $\{h_j\}_{j=1}^n$ in $\diffloc{1+\alpha}(M)$ with $h_j:\ball{x_{j-1}}{R}\to h_j(\ball{x_{j-1}}{R})$  satisfying $h_j(x_{j-1})=x_j$, $\|Dh_j\|<\lambda$, $\|h_j\|_{\cC^{1+\alpha}}<C$, and every pair of measurable sets $S_1,S_2\subseteq \ball{x_0}{R}$ of positive Lebesgue measure,
    \[L^{-1}\frac{\leb(S_1)}{\leb(S_2)}<\frac{\leb(h^n(S_1))}{\leb(h^n(S_2))}<L\frac{\leb(S_1)}{\leb(S_2)}.\]
    \end{lemma}
    Recall that $h^n:=h_n\circ \cdots \circ h_1$.
    \begin{proof}
    Since there is an upper bound for the $\cC^{1+\alpha}$ norm of the derivative of the $\exp$ function on the balls of radius $R_1$ on the whole manifold $M$, by replacing $h_j$ with $\exp_{x_j}\circ h_j\circ \exp_{x_{j-1}}^{-1}:\Eucball{R}\to \Eucball{R}$, one can assume that the maps are defined between open sets of $\R^d$. Now, it is enough to show that there exists $L_1>1$ such that for any measurable set $S\subseteq \Eucball{R}$, 
    \begin{equation}\label{eq:vol-det}
    L_1^{-1}|\det D_0h^n|{\leb(S)} \leq {\leb(h^n(S))} \leq L_1|\det D_0h^n|{\leb(S)}.   
    \end{equation}
    Since the sequence has bounded $\cC^{1+\alpha}$ norm, the maps $z\mapsto \log|\det D_zh_j|$ are $\alpha$-H\"{o}lder on $\Eucball{R}$ with some uniform constant, that is, there exists $L'>0$ (independent of $j$) such that for any $j\geq 1$ and any pair $z,z'\in \Eucball{R}$, 
    \[\Big|\log \big|\det D_{z}h_j\big|-\log \big|\det D_{z'}h_j\big|\Big|<L'|z-z'|^\alpha.\]
    For $x\in \Eucball{R}$ and $j\leq n$, denote $x_j:=h^j(x)$. From the contraction property, 
    \[|x_j|\leq \Big(\sup\limits_{z\in \Eucball{R}}\|D_zh^j\|\Big)|x|\leq \lambda^j|x|\leq \lambda^j R.\]
    Therefore,
    \begin{align*}
    \Big|\log \frac{|\det D_{x}h^n |}{|\det D_{0}h^n |}\Big|
    &=\sum\limits_{j=0}^{n-1}\Big| \log|\det D_{x_j}h_{j+1}|-|\det D_{0}h_{j+1}|\Big|\\
    &< L'\sum\limits_{j=0}^{n-1}|x_j|^\alpha \leq L'R^\alpha \sum\limits_{j=0}^{n-1}  \lambda^{j\alpha}.
    \end{align*}
    Now, since $\uppercont<1$, $L_1:=\exp\big(L'R^\alpha\sum\limits_{j=0}^\infty \lambda^{j\alpha}\big)<\infty$ and so
    \begin{equation}\label{eq:control-det}
    L^{-1}_1|\det D_0h^n|\leq |\det D_xh^n|\leq L_1|\det D_0h^n|.   
    \end{equation}
    By the change of variable formula,
    \begin{align*}
        \int_S|\det D_xh^n|\mathrm{d}\leb(x)\leq  L_1|\det D_0h^n|\leb(S).
    \end{align*}
    The proof of the other inequality in (\ref{eq:vol-det}) is similar.
    \end{proof}
    \subsection{Infiltrated  quasi-conformality}
    In this subsection, we show that under hypotheses (H0)-(H3), quasi-conformality of the derivatives at the origin leads to the quasi-conformality in a neighbourhood. Informally, the idea is that the contracting assumption forces the derivatives of long blocks in the nearby points to imitate the behaviour of derivatives at the origin. The results of this subsection will not be used in other parts of the paper and are included here for their own interest. 
  
    \begin{theorem}\label{thm:st-control-shape}
     Let  $x\in M$, $\{f_i\}_{i=1}^\infty$ be a sequence in $\mathrm{Diff}^{1+\alpha}(M)$ with bounded $\cC^{1+\alpha}$ norm and  quasi-conformal expanding at $x$. Then, there exist $R>0$ and $\theta>1$ such that for any $n\in \N$ and any ball $\ball{y}{r}\subseteq \ball{f^n(x)}{R}$,
    \[ \ball{f^{-n}y}{r_n}\subseteq f^{-n} (\ball{y}{r})\subseteq    \ball{f^{-n}(y)}{\theta r_n},\]
      where $r_n= r\theta^{-\frac{1}{2}}  |\det D_{x}f^{-n}|^{1/d} $.
    \end{theorem}
    To prove Theorem \ref{thm:st-control-shape}, we first prove the following proposition. 
  
    \begin{proposition}\label{prop:infiltration}
    Let $R,C>0$, $\kappa\geq 1 >\uppercont> \lowercont>0$ and  $\alpha\in (0,1)$. Then, there exist $\xi_0>0$, $\overline{\kappa}>1$ such that any sequence $\{h_j\}_{j=1}^\infty$ of maps satisfying hypotheses (H0)-(H3) is $\overline{\kappa}$-conformal at every point of $\Eucball{\xi_0}$. 
    \end{proposition}
    \begin{proof}
    The proposition follows from a refinement of the proof of Lemma \ref{lem:pseudo-orbit}. For $\alpha<1$, one should take blocks of compositions and repeat the claims of the proof of Lemma \ref{lem:pseudo-orbit}. To avoid repeating the arguments, here we give a proof for $\cC^{1+\mathrm{Lip}}$ regularity, that is, for $\alpha=1$. For $x\in \Eucball{R}$, denote $x_i:=h^i(x)$, $D_i:=D_0h_i$, $\delta_i:=|\det D_i|^{1/d}$,  $\oD_i:=D_{x_{i-1}}h_i$ and $\cE_i:=D_i-\oD_i$. Observe that  for $n\in \N$,
    \[D_0h^n-D_xh^n=D_n\cdots D_1-\oD_n\cdots \oD_1=\sum\limits_{i=1}^{n} D_n\cdots D_{i+1} \cE_i \oD_{i-1}\cdots \oD_1.\]
    Hence,
    \begin{equation}\label{eq:difference-of-derivatives-1}
    \|D_xh^n-D_0h^n\|\leq \sum\limits_{i=1}^{n} \|D_n\cdots D_{i+1}\| \cdot \|\cE_i\| \cdot  \|\oD_{i-1}\cdots \oD_1\|. \end{equation}
    Since $D_n\cdots D_{i+1}$ is $\kappa^2$-conformal, in view of (C2)$^\prime$, one gets that $\|D_n\cdots D_{i+1}\|\leq \ozetatwo \delta_{n}\cdots \delta_{i+1}$. On the other hand, (H1) implies $\|\oD_{i-1}\cdots \oD_1\|\leq \uppercont^{\,i-1}$.
    Now, for $\theta>1$ given by Lemma \ref{lem:iteration-balls-Rd}, $|x_{i-1}|\leq \theta\delta_{i-1}\cdots \delta_1 |x|$ and so, by (H3), \[\|\cE_i\|= |D_{x_{i-1}}h_i-D_0h_i|\leq C \theta |x_{i-1}| \leq C\theta R \delta_{i-1}\cdots\delta_1 .\]
    Using (\ref{eq:difference-of-derivatives-1}), 
    \begin{align*}
   \|D_xh^n-D_0h^n\|&\leq \sum\limits_{i=1}^{n} C R \ozetatwo \delta_n\cdots\delta_{i+1} \delta_{i-1}\cdots \delta_1 \uppercont^{\,i-1}\\
   &\leq CR\ozetatwo \lowercont^{-1} \Big(\sum\limits_{i=1}^n {\uppercont^{\,i-1}}\Big) \delta_n \cdots \delta_1 .
    \end{align*}
    Now, by the convergence of the series $\sum_{i=0}^\infty \uppercont^{\,i-1}$, there exists $C_1>0$ such that for any $n\geq 1$ and $x\in \Eucball{R}$, 
    \begin{equation}\label{eq:close-derivatives}
    \|D_xh^n\|\leq \|D_0h^n\|+\|D_xh^n-D_0h^n\|<C_1 |\det D_0h^n|^{1/d}.
    \end{equation}
    By (\ref{eq:control-det}), one obtains that $\|D_xh^n\|\leq C_1L_1^{1/d}|\det D_xh^n|^{1/d}$. Thus, there exists $\overline{\kappa}=\overline{\kappa}(d,C_1L^{1/d}_1)>0$ such that $D_xh^n$ is $\overline{\kappa}$-conformal, as claimed.  
    \end{proof}
      
    \begin{proof}[Proof of Theorem \ref{thm:st-control-shape}]
    The proof is similar to the one of Theorem \ref{thm:control-shape-M}. Take $x_i,\upperexp,\lowerexp,R$ as in the proof of that theorem. Define $\{y_i\}_{i=0}^n$ by $y_n:=y$ and $y_{i-1}:=f_i^{-1}(y_i)$. By uniform expansion, one obtains that for any $j< n$, \[f^{-1}_{j+1}\circ \cdots \circ f^{-1}_n(\ball{y_n}{r})\subseteq \ball{y_j}{r}\cap \ball{x_j}{R}.\] Now, define $\hat{f}_i:=(\exp_{y_i}^{-1}\circ f_i\circ \exp_{y_{i-1}})|_{\Eucball{r}}$ and $\hat{h}_i:=\hat{f}^{-1}_{n+1-i}|_{\Eucball{r}}$. So, the conclusion follows from Lemma \ref{lem:iteration-balls-Rd} for this sequence. Indeed, condition (H3) is guaranteed by Proposition \ref{prop:infiltration}. 
    \end{proof}
    
    \begin{remark}
    By Proposition \ref{prop:infiltration}, Lemma \ref{lem:distortion-part1} and the results in the theory of quasi-conformal and quasi-symmetric maps, one can give another proof for Theorem \ref{thm:st-control-shape}. In fact, if a map is $\kappa$-conformal on its domain, then there exists a bound for the ratio between outer and inner radii of the image of a ball (see \cite[Section 4]{Heinonen_Koskela_1998} and \cite{Vaisala}). Then, the inner and outer radii can be estimated by means of estimating the volume and the bounded distortion lemma (Lemma \ref{lem:distortion-part1}). 
    \end{remark}

    \section{Quasi-conformal blenders}\label{sec:blender}
    This section is devoted to a new mechanism/phenomenon that we call quasi-conformal blender.  For pseudo-semigroup actions, the quasi-conformal blender guarantees the existence of quasi-conformal expanding orbit-branches at every point in some region, leading to the stable local ergodicity. Here, we present the proof of Theorem \ref{thm:blender} and its variants using the results of the previous sections.
   
    Throughout the section, $M$ is a boundaryless, not necessarily compact, smooth Riemannian manifold of dimension $d$. 
    
 \begin{definition}[$\rho$-ergodic]
    Let $V\subseteq M$ be an open set with compact closure. For $\rho>0$ and $\cF\subseteq \diffloc{1}(M)$, we say $\IFS(\cF\loc{V})$ is  \emph{$\rho$-ergodic} (w.r.t.  Leb.), if every measurable $\cF\loc{V}$-invariant set of positive measure in $V$ contains a ball of radius $\rho$, up to a set of zero Lebesgue measure.
    \end{definition}

    This definition is equivalent to saying that the set of density points of every measurable $\cF\loc{V}$-invariant subset of $V$ either is empty or contains a ball of radius $\rho$.   
    Recall that $x\in M$ is a \emph{(Lebesgue) density point} of a measurable set $S\subseteq M$ if
    \begin{equation}\label{eq:density-point}
    \lim\limits_{r\to 0}\frac{\leb\big(S\cap \ball{x}{r}\big)}{\leb\big(\ball{x}{r}\big)}=1.
    \end{equation}

    \subsection{From quasi-conformal expansion to local ergodicity}
    In this subsection, we state a technical lemma about the quasi-conformal expanding sequences. It will be used in both local and global settings for proving local ergodicity.

    \begin{lemma}[Local ergodicity] \label{lem:local-erg} 
    Let
    $\{f_i\}_{i=1}^\infty$ be a sequence of expanding diffeomorphism in $\diff_{\rm loc}^{1+\alpha}(M)$ with bounded $\cC^{1+\alpha}$ norm. 
    Let also $x_0\in M$ and $\rho>0$ be such that for any $n\in \N$, $B(f^{n-1}(x_0),\rho)\subset \mathrm{Dom}(f_n)$. If the sequence is quasi-conformal at $x_0$ and $\{f^i(x_0)\}_{i\geq 1}$ is bounded in $M$, then 
    \begin{itemize}
    \item[(a)] for any open set $U\subseteq \ball{x_0}{\rho}$ containing $x_0$, there exists $n\in \N$ with $\ball{f^n(x_0)}{\rho}\subseteq f^n(U)$,
    \item[(b)] for any measurable set $S\subseteq \ball{x_0}{\rho}$ with density point at $x_0$, there exists $n\in \N$ such that $\bigcup_{i\in \N}f^i(S)$ contains $\ball{f^n(x_0)}{\rho}$, up to a set of zero Lebesgue measure. 
    \end{itemize}
    \end{lemma}
    \begin{proof} Let $\kappa,\eta>1,C>0$ be such that the sequence $\{f_i\}_{i=1}^\infty$ is $\eta$-expanding,  $\kappa$-conformal at $x_0$ and  $\sup_j\|Df_j\|_{\cC^{1+\alpha}}<C$.
    Also, denote $x_n:=f^n(x_0)$. 
    
    \medskip
    
    \noindent
    \textit{Proof of ${\mathrm{(a)}}$.} Fix open set $U\subseteq B(x_0,\rho)$ containing $x_0$. For every $i\geq 1$, let $s_i$ be the largest positive number in $(0,\rho]$ satisfying $\ball{x_i}{s_i}\subseteq f^i(U)$. Since $m(D{f}_i|_{\ball{x_{i-1}}{\rho}})>\eta$, if for some $i\geq 1$, $s_i<\rho \eta^{-1}$, then $s_{i+1}>\eta s_i$, and if $s_{i}\geq \rho\eta^{-1}$, then $s_{i+1}=\rho$. Hence, there is $n\geq 0$ with $s_n=\rho$. This finishes the proof of part (a). Note that for this part of the lemma, the sequence $\{f_i\}$ only needs to be expanding and $\cC^1$-regular. 
    
    \medskip
    \noindent
    \emph{Proof of $\mathrm{(b)}$}. Let  $S\subseteq B(x_0,\rho)$  be a measurable set with a density point at $x_0$. By applying Theorem \ref{thm:control-shape-M} to the sequence $\{f_i\}_i$, one can find $\xi_0>0$ and $\theta>1$
    such that for each $j\in \N$, there is $r_j>0$ with    
    \begin{equation}\label{eq:good-n}
    f^j(\ball{x_0}{r_j})\subseteq \ball{x_j}{\xi_0}  \subseteq f^j(\ball{x_0}{\theta r_j}),
    \end{equation} 
    and $\lim\limits_{j\to\infty}r_j=0$. 
    For the rest of the proof, fix $\xi_0$ and assume that $\xi_0<\rho$.  Denote $\hat{S}:=\bigcup_{i\geq 0}f^i(S)$. Since $x_0$ is a density point of $S$, for any $\epsilon>0$, there exists $j_0\in\N$ such that whenever $j>j_0$,
    \[\frac{\leb \big( \ball{x_0}{\theta r_j}\setminus \hat{S}\big)}{\leb \big(\ball{x_0}{\theta r_j}\big)}<{\epsilon}.\] 
    There exists $\sigma>0$ such that for any $r,r'\in(0,\rho)$, 
    \[\frac{\leb(\ball{x_0}{r})}{\leb(\ball{x_0}{r'})}\leq \sigma (\frac{r}{r'})^d.\]
    Denote $f^{-j}:=f_1^{-1}\circ \cdots \circ f_j^{-1}$. Then, by (\ref{eq:good-n}),
    \begin{eqnarray*}
    \frac{\leb\big( f^{-j}\big(\ball{x_j}{\xi_0}\setminus \hat{S})\big)}{\leb \big(f^{-j}(\ball{x_j}{\xi_0})\big)}&\leq  &\frac{\leb\big( \ball{x_0}{\theta r_j}\setminus \hat{S}\big)}{\leb \big( \ball{x_0}{\theta r_j}\big)}\cdot\frac{\leb \big( \ball{x_0}{\theta r_j}\big)}{\leb \big(f^{-j}(\ball{x_j}{\xi_0})\big)}\\
    &\leq &\frac{\leb\big(\ball{x_0}{\theta r_j}\setminus \hat{S}\big)}{\leb \big( \ball{x_0}{\theta r_j}\big)}\cdot\frac{\leb \big( \ball{x_0}{\theta r_j}\big)}{\leb ( \ball{x_0}{r_j}\big)}\\
    & < & {\epsilon}\sigma\theta^d.
    \end{eqnarray*}
    It follows from Lemma \ref{lem:distortion-part1} that for some $L> 1$, 
    \[\frac{\leb\big(\ball{x_j}{\xi_0}\setminus \hat{S}\big)}{\leb\big(\ball{x_j}{\xi_0}\big)}<{L \frac{\leb\big(f^{-j}(\ball{x_j}{\xi_0}\setminus\hat{S})\big)}{\leb\big(f^{-j}(\ball{x_j}{\xi_0})\big)}<} {\epsilon}\sigma L\theta^d.\]
    Thus, for any $j> j_0$,
    \begin{equation}\label{eq:high-density}
        \frac{\leb\big(\hat{S}\cap \ball{x_j}{\xi_0}\big)}{\leb\big(\ball{x_j}{\xi_0}\big)}> 1-{\epsilon}\sigma L\theta^d.
    \end{equation}

    Now, since $\epsilon$ was arbitrary, by (\ref{eq:high-density}), the density of $\hat{S}$ in $\ball{x_j}{\xi_0}$ tends to $1$ (as $j\to\infty$). Then, for each accumulation point $y_0$ of the bounded sequence $\{x_i\}_{i=1}^\infty$,  $\hat{S}$ contains $\ball{y_0}{\xi_0}$ up to a set of zero Lebesgue measure. {Next, take a sufficiently large $i$ such that $x_i\in \ball{y_0}{\xi_0}$. This implies that $\hat{S}$ contains an open neighbourhood $U$ of $x_i$, up to a set of zero Lebesgue measure. Then, by part (a), there exists $n>i$ such that $\ball{x_n}{\rho}\subseteq f^{n-i}(U)$.   Finally, since diffeomorphisms maps sets of zero Lebesgue measure to sets of zero Lebesgue measure, $\hat{S}$ contains $\ball{x_n}{\rho}$, up to a set of zero Lebesgue measure.} 
    \end{proof}
    
    Lemma \ref{lem:local-erg} has the following global consequence which can be seen as a generalization of Theorem \ref{thm:sullivan}.
    
    \begin{theorem}\label{thm:gen-sullivan}  
    Let $M$ be a closed manifold and $\mathcal{F}\subseteq \mathrm{Diff}^{1+\alpha}(M)$ be finite. Suppose that there exist $\eta,\kappa>1$ such that $\ifs{F}$ has a $\kappa$-conformal $\eta$-expanding orbit-branch at every point. Then,   $\ifs{F}$ is $\rho$-ergodic for some $\rho>0$.
    In particular, the  action of every group $G\subseteq \diff^{1}(M)$ containing $\cF$ is ergodic, provided that it is minimal.
    
    \end{theorem}
    \begin{proof}
    Since $\cF$ is finite, there exist $\rho>0$ and $\eta'\in (1,\eta)$ such that whenever $m(D_xf)>\eta$, for some $x\in M$ and $f\in \cF$, then $m(Df|_{\ball{x}{\rho}})>\eta'$. Indeed, if $C:=\max_{f\in \cF}\|f\|_{\cC^{1+\alpha}}$, then by \eqref{eq:diff-conorm}, for any $x,y\in M$,
    \[|m(D_xf)-m(D_yf)|\leq Cd(x,y)^\alpha.\]
    where $d(.,.)$ denotes the distance on $M$. So, $m(D_xf)>\eta$ implies that $m(D_yf)>\eta'$, provided that $d(x,y)^\alpha<\rho:=C^{-1}(\eta-\eta')$.
    
    For $f\in \cF$, denote $U_f:=\{x\in M{\gap}m(D_xf)>\eta'\}$ and $\hat{f}:=f|_{U_f}$. Also, let $\hat{\cF}:=\{\hat{f}{\gap}f\in \cF\}$. Consider a measurable $\cF$-invariant subset $S$ of positive measure and pick $x$ to be a density point of $S$. Then, the $\kappa$-conformal $\eta$-expanding orbit-branch of $\IFS(\hat{\cF})$ at $x$, provides a sequence of maps satisfying the assumptions of Lemma \ref{lem:local-erg} and the conclusion follows from this lemma. 
    
    For the second part, denote the set of density points of $S$ by $S^\bullet$. Since $\cF\subseteq G$,  it follows from the first part that $S^\bullet$ contains an open ball $B$. We claim that $S^\bullet=M$. Indeed, 
    $\leb(B\setminus S)=0$ implies that  for every $g\in G$,  $\leb(g(B)\setminus g(S))=0$. Then, by the invariance of $S$, $\leb(g(B)\setminus S)=0$ and in particular, $g(B)\subseteq S^\bullet$. On the other hand, by the minimality assumption, $\bigcup_{g\in G}g(B)=M$. This proves the claim and finishes the proof of the theorem. 
    \end{proof}

    \subsection{Proof of Theorem \ref{thm:blender}}
    We will prove the following theorem which in particular implies Theorem \ref{thm:blender}.

    \begin{theorem}[Quasi-conformal blender] 
    \label{thm:blender-full}
    Let $\cF \subseteq \diff_{\rm loc}^{1+}(M)$. 
    Let $\cW \subseteq \cE(M)$ be an open set with compact closure and $V:=\pi(\cW)$. Assume that  
    for any ${\bf w}\in \overline{\cW}$, there exists $f\in \cF$ satisfying
        \begin{itemize}
            \item[(i)] $\ef f({\bf w})\in\cW$,
            \item[(ii)]  $m(D_xf)>1$, where $x=\pi({\bf w})$. 
        \end{itemize}
    Then, there exist real numbers $\rho>0$ and $\kappa>1$ such that for every $\tilde{\cF}\subseteq\diff_{\rm loc}^{1+}(M)$ sufficiently close to $\cF$ in the $\cC^1$ topology, 
        \begin{itemize}
            \item[{(a)}] for any $x\in V$, $\IFS(\tilde{\cF}\loc{V})$ has an orbit-branch which is $\kappa$-conformal at $x$,
            \item[{(b)}] $\IFS(\tilde{\cF}\loc{V})$ is $\rho$-ergodic.
        \end{itemize}
    In addition, if $M$ is compact, $\IFS(\cF)$ and $\IFS(\cF^{-1})$ are minimal, 
    then $\IFS(\cF)$ is $\cC^1$-stably ergodic in $\diffloc{1+}(M)$ and $\IFS(\cF^{-1})$ is $\cC^1$-robustly minimal.
    \end{theorem}
    
    \begin{figure}[t]
    \begin{minipage}{.40\textwidth}
    
    \centering
    \includegraphics[width=.9\textwidth]{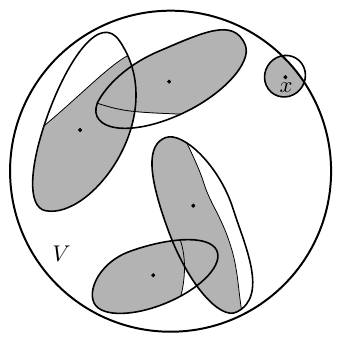}
    Images of a neighbourhood of $x$
    \end{minipage}
    \begin{minipage}{.1\textwidth}
    $~$
    \end{minipage}
    \begin{minipage}{.4\textwidth}
    \centering
    \vspace{.5cm}
    \includegraphics[width=.75\textwidth]{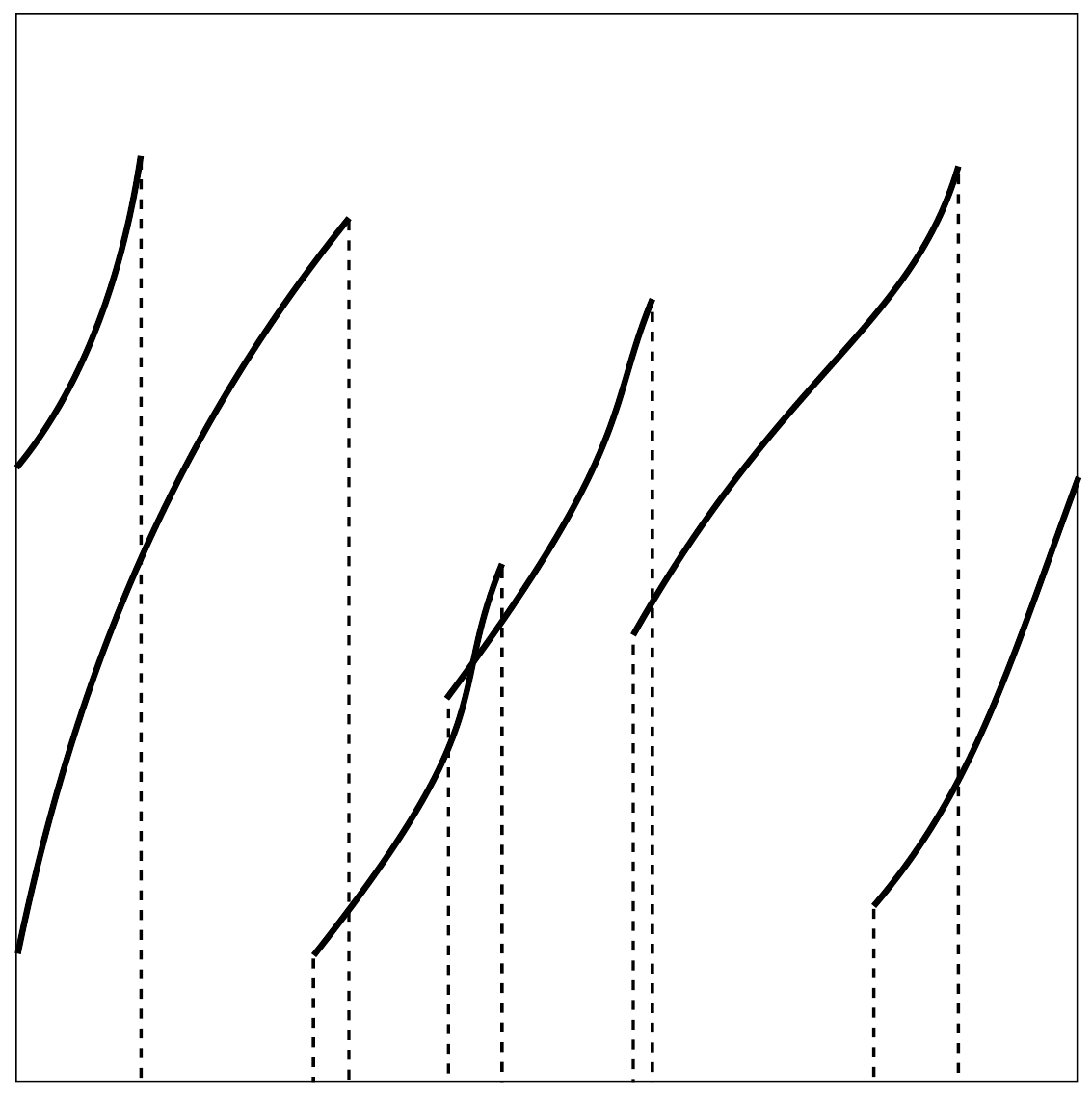}
    
    \vspace{.4cm}
    
    One-dimensional localized maps
    \end{minipage}
    \caption{Families of expanding maps satisfying the covering condition in Theorem \ref{thm:blender}.}
    \label{fig:blender}
    \end{figure}
    
    \begin{proof}[Proof of Theorem \ref{thm:blender}]
    It follows from (\ref{eq:covering}) that for every ${\bf w}\in \cW$, there exists $f\in \cF$ with $\ef f ({\bf w})\in \cW$. Now, since every element of $\cF$ is expanding, one has $m(D_xf)>1$ for $x=\pi({\bf w})$. So, Theorem \ref{thm:blender} follows from part (b) of Theorem \ref{thm:blender-full}.  
    \end{proof}
    \begin{remark}
    If a family $\cF$ and a set $\cW$ satisfy the assumptions of Theorem \ref{thm:blender-full}, then for ${\bf w}\in\overline{\cW}$, one can get a subfamily $\cF_{\bf w}\subseteq \cF$ consisting of all elements satisfying (i)-(ii). Then, for any element $f\in \cF_{\bf w}$, restrict its domain to a small neighbourhood of ${\bf w}$ such that the restricted map is expanding. Let $\cF'$ be the family consisting of all these restricted diffeomorphisms. It is clear that $(\cF',V)$ is a quasi-conformal blender. In other words, one can deduce part (b) in Theorem \ref{thm:blender-full} from  Theorem      \ref{thm:blender}, and vice-versa.
    \end{remark}

    \begin{remark}
    As a matter of fact, the $\cC^1$ stability in Theorem \ref{thm:blender-full} and its consequences in this paper are valid in a substantially stronger form that allows to perturb the family at each step of iterations. We do not discuss the details in this paper. Cf. \cite{homburg_nassiri_2014}, where this notion of {\it strong stability} has been introduced.
    \end{remark}
     
    \begin{proof}[Proof of Theorem \ref{thm:blender-full}]
    Since the assumptions (i)-(ii) are stable under small perturbation of $f$, in the $\cC^1$ topology, and ${\bf w}\in \cE(M)$, one can deduce that 
    \begin{itemize}
         \item For every ${\bf w}\in \overline{\cW}$ there exist $\epsilon({\bf w})>0$ and $f_{\bf w}\in \cF$ such that if $\cB_{\bf w}$ denotes the open ball of radius $\epsilon({\bf w})$ with center ${\bf w}\in\cE(M)$, then for any $\tilde{f}$ sufficiently close to $f_{\bf w}$, in the $\cC^1$ topology, $\pi(\cB_{\bf w})\subseteq {\rm Dom}(\tilde{f})$ and for any ${\bf w}'\in \cB_{\bf w}$,
         \begin{itemize}
             \item[(1)] $\ef \tilde{f}({\bf w}')\in \cW$,
             \item[(2)] $m(D_{x'}\tilde{f})>1$, where $\pi({\bf w}')=x'$.
         \end{itemize}
    \end{itemize}
    
    \medskip
    \noindent
    \textit{Proof of $\rm{(a)}$.}
    Let $\cB'_{\bf w}$ be the open ball of radius  $\frac{1}{2}\epsilon({\bf w})$ with center ${\bf w}$. Then, by the compactness of $\overline{\cW}$, there is a finite subset $\{{\bf w}_1,\ldots,{\bf w}_k\}$ of $\overline{\cW}$ such that $\overline{\cW}\subseteq \bigcup_{i}\cB'_{{\bf w}_i}$. By Corollary \ref{cor:st-QC}, this implies the existence of $\kappa>1$ such that for every $\tilde{\cF}$  sufficiently close to $\cF$ in the $\cC^1$ topology, $\IFS(\tilde{\cF}\loc{V})$ has a $\kappa$-conformal orbit-branch at every point of $V$. So, the proof of part (a) is finished. 
    
    \medskip
    \noindent
    \textit{Proof of $\rm{(b)}$.} Let $\rho>0$ be smaller than the Lebesgue number of the open covering $\bigcup_{i}\cB'_{{\bf w}_i}$ for $\overline{\cW}$. Then, there is $\eta>1$ such that for any ${\bf w}\in \overline{\cW}$, there exists $1\leq i \leq k$ such that for any $\tilde{f}$ sufficiently close to $f_{{\bf w}_i}$, in the $\cC^1$ topology, $\ball{x}{\rho}\subseteq \mathrm{Dom}(\tilde{f})$, $\tilde{f}(\ball{x}{\rho})\subseteq V$, and   $m(D\tilde{f}|_{\ball{x}{\rho}})>\eta$, where $x=\pi({\bf w})$.
    
    Denote $\tilde{\cF}_0:=\{\tilde{f}_{{\bf w}_1},\ldots,\tilde{f}_{{\bf w}_k}\}$. Let $\alpha>0$ be such that $\tilde{\cF}_0\subseteq \diffloc{1+\alpha}(M)$. Consider a measurable $\tilde{\cF_0}\loc{V}$-invariant set $S$ of positive measure and let $x_0$ be a density point of $S$. Family $\cF_0$ satisfies the assumptions (i)-(ii) of the theorem, so it follows from part (a) that $\IFS(\tilde{\cF}_0|_V)$ has a $\kappa$-conformal orbit-branch $\{x_i\}_{i=0}^\infty$ at $x_0$. Let $\{\tilde{f}_i\}_{i=1}^\infty$ be the sequence of maps providing this orbit-branch, namely $f_i(x_{i-1})=x_i$ for every $i\in \N$. Next, our aim is to apply Lemma \ref{lem:local-erg} to this sequence. The problem is that $x_0$ may be close to the boundary of $V$ and $\ball{x_0}{\rho}\not\subseteq \tilde{f}_1\loc{V}$. To avoid this challenge, we remove the first term of the sequences and consider $\{\tilde{f}_i\}_{i=2}^\infty$ and $\{x_i\}_{i=1}^\infty$, which by means of above arguments satisfy the assumptions of Lemma \ref{lem:local-erg}. Note that $x_1=\tilde{f}_1(x_0)$ is also a density point of the invariant set $S$. By part (b) of  Lemma \ref{lem:local-erg}, there is $\xi_0>0$ and an accumulation point  $y_0$ of $\{x_i\}_{i=1}^\infty$ such that $S$ contains $\ball{y_0}{\xi_0}$, up to a set of zero Lebesgue measure. Since  for any $i\geq 1$, $\ball{x_i}{\rho}\subseteq V$ and so $\ball{y_0}{\rho}\subseteq V$. Similarly, an expanding orbit-branch  $\{y_i\}_{i=0}^\infty$ of $\IFS(\tilde{\cF}\loc{V})$ at $y_0$ can be provided in such a way that for any $i\geq 0$, $\ball{y_i}{\rho}\subseteq V$. Using part (a) of Lemma \ref{lem:local-erg} for open set $U=\ball{y_0}{\xi_0}$, one can conclude that the set of density points of $S$ contains $\ball{y_n}{\rho}$ for some $n\geq 0$. This finishes the proof of part (b). 

    \medskip
    
    Next, we move on to the proof of global results. We assume that $M$ is compact, $\IFS(\cF)$ and $\IFS(\cF^{-1})$ are minimal on $M$.  By the minimality, $\langle {\cF}\rangle^+(\ball{x}{\rho'})=\langle {\cF}^{-1}\rangle^+(\ball{x}{\rho'})=M$  where $\rho'<\frac{1}{2}\rho$. Suppose that $\bigcup_{x\in X}\ball{x}{\rho'}=M$ for some finite set $X\in M$. On the other hand, by the compactness of $M$, there exists a finite set $\cF_1\subseteq \langle {\cF}\rangle^+$ with $\cF_1(\ball{x}{\rho'})=\cF_1^{-1}(\ball{x}{\rho'})=M$, for any $x\in X$. Consider small perturbation $\tilde{\cF}$ of $\cF$ and denote the elements of $\tilde{\cF}$ corresponding to the family $\cF_1$ by $\tilde{\cF}_1$. If the perturbation is sufficiently small, then for every ball $B$ of radius $\rho$,
    \begin{equation}\label{eq:rho-min}
    \tilde{\cF}_1(B)=\tilde{\cF}_1^{-1}(B)=M.
     \end{equation}
    
    \noindent
    \textit{Proof of stable ergodicity.}  Let $S$  be a  measurable $\tilde{\cF}$-invariant subset of $M$ with positive measure.  Consider an arbitrary ball $B_0\subseteq V$ of radius $\rho$. By (\ref{eq:rho-min}), $S\subseteq \tilde{\cF}_1^{-1}(B_0)=M$ and so  $\leb(S\cap B_0)>0$. Then, by $\rho$-ergodicity of $\IFS(\tilde{\cF}\loc{V})$, the set of density points of $S$ contains some ball $B_1\subseteq M$ of radius $\rho$. Finally, by (\ref{eq:rho-min}), $\tilde{\cF}_1(B_1)=M$ and this implies that $S$ contains $M$, up to a set of zero Lebesgue measure.
    
    \medskip
    \noindent
    \textit{Proof of robust minimality.} The proof of robust minimality is similar to the one of stable ergodicity.  Consider an open set $U\subseteq M$, and ball $B_0$ as above. Let $\tilde{\cF}\subseteq \diffloc{1}(M)$ be
     sufficiently close to $\cF$ in the $\cC^1$ topology,
    satisfying (\ref{eq:rho-min}). Again, (\ref{eq:rho-min}) implies that $\hat{U}:=\langle\tilde{\cF}\rangle^+(U)$ intersects $B_0\subseteq V$. Consider $z_0\in \hat{U}\cap B_0$. Then, similar to the arguments of part (b), one can find an expanding orbit-branch of $\IFS(\tilde{\cF}\loc{V})$ at $z_0$ and use part (a) of Lemma \ref{lem:local-erg} to deduce that $\hat{U}$ contains a ball $B_2\subseteq V$ of radius $\rho$.  Finally, again by (\ref{eq:rho-min}), $M= \tilde{\cF}_1(B_2)\subseteq \hat{U}$. Since open set $U$ was arbitrarily chosen, it follows that the orbit of every point under $\IFS(\tilde{\cF}^{-1})$ is dense in $M$. This means $\IFS(\tilde{\cF}^{-1})$ is minimal. 
    \end{proof}
    
    \subsection{Contracting quasi-conformal blenders} 
    In this subsection, we present a variant of Theorem \ref{thm:blender} in the setting of contracting maps. This is useful to get local ergodicity and local minimality simultaneously.
   
    \begin{theorem}
    \label{thm:blender-cont}
    Let $\cG\subseteq \diffloc{1+}(M)$ and $\cW\subseteq \cE(M)$ be a bounded open set with $V:=\pi(\cW)$. Let $U\subseteq M$ be an open set with compact closure containing $\pi(\overline{\cW})$. Also, assume that 
    \begin{itemize}
        \item[(i)]   $\overline{\cW}\subseteq \bigcup_{g\in \cG}\ef g (\cW)$,
        \item[(ii)] for every $g\in\cG$, $\overline{U}\subseteq \mathrm{Dom}(g)$ and $g(\overline{U})\subseteq U$. 
        \item[(iii)] every $g\in\cG$ is uniformly contracting on $U$.
        \end{itemize}
        Then,  for every $\tilde{\cG}\subseteq\diffloc{1+}(M)$ 
         sufficiently close to $\cG$ in the $\cC^1$ topology,
        any measurable $\tilde{\cG}\loc{U}$-invariant set of positive measure in $U$, has full measure in $V$. Also, for every $x\in U$, $\langle\cG\loc{U}\rangle^+(x)$ is dense in $V$.
    \end{theorem}
    \begin{proof}
    By the compactness of $\overline{\cW}$ and the openness of $\cW$, one can replace $\cG$ with a finite subfamily satisfying the assumptions. So, we may assume that $\cG\subseteq \diffloc{1+\alpha}(M)$, for some $\alpha>0$, is finite. 
    
    Let $\Lambda:=\cap_n\cG^n(\overline{U})$ be the Hutchinson attractor of $\IFS(\cG|_{U})$.
    It follows from \cite{Hutchinson} that for any $x\in U$, ${\gen{\cG}^+(x)}$ is dense in $\Lambda$ (cf. \cite[Theorem 4.2]{homburg_nassiri_2014}). On the other hand, for any $n\in \N$,
    \[V\subseteq \cG(V) \subseteq \cG^n(V)\subseteq \cG^n(\overline{U}),\]
    and thus $\overline{V}\subseteq \Lambda$. This, in particular, implies that for every $x\in U$, $\langle\cG\loc{U}\rangle^+(x)$ is dense in $V$ and consequently, every $\cG\loc{U}$-invariant subset of positive measure in $U$, intersects $V$ in a set of positive measure. 
    
    Note that the family $\cG^{-1}\loc{U}$ is uniformly expanding on $V$. On the other hand, in view of assumption (i), $\cG^{-1}$ satisfies the covering property (\ref{eq:covering}) for $\cW$. Now, let $S$ be a measurable $\cG\loc{U}$-invariant set with $\leb(S)>0$. It follows from above that $\leb(S\cap V)>0$. Suppose that $\leb(S\cap V)<\leb(V)$.  Then, $S':=V\setminus S$ is  $\cG^{-1}\loc{V}$-invariant and $\leb(S')>0$. By applying Theorem \ref{thm:blender-full} to the family $\cG^{-1}\loc{V}$, one gets that the set of density points of $S'$ contains an open ball $B$. Let $x$ be a density point of $S$. By the above arguments, some element of $\langle\cG\loc{U}\rangle^+$ maps $x$ to $B$. This contradicts the invariance of $S$ under $\cG\loc{U}$ and shows that $S$ must have full measure in $V$.
    
    Finally, note that after a small perturbation of the generators of $\cG$, all the assumptions of the theorem hold. That is, the conclusion follows for the family $\tilde{\cG}$ sufficiently close to $\cG$, following similar arguments above. 
    \end{proof}
    
    \begin{remark}
    It follows from the proof of Theorem \ref{thm:blender-cont} that one can replace the assumption (i) with $\overline{V}\subseteq \bigcup_{g\in \cG}g(V)$ and deduce the density of $\langle\cG\rangle^+(x)$ in $V$ for any $x\in V$.
    \end{remark}
    
    \subsection{Other statements on quasi-conformal blenders}
    In this subsection, we state an analogue of Theorem \ref{thm:blender}.  
    First note that if $V\subseteq M$ is completely inside an open chart of $M$, one can use the coordinating map to translate the problem to an open subset of the Euclidean space. In this case, one can get the following statement in view of Theorem \ref{thm:blender-full}.  
    
    \begin{theorem}\label{thm:blender-Rd} 
        Let $\mathcal{F}\subseteq \diff^{1+}_{\rm loc}(\R^d)$ be a family of orientation-preserving maps. 
        Let  $V\subseteq \R^d$ be a bounded open set and  $\cU \subseteq \sldr$ be a bounded open neighbourhood of the identity. Assume that for any $x\in \overline{V}$, there exists $\cF_x\subseteq \cF$ such that
        \begin{itemize}
            \item[(i)]  $f(x)\in V$, for $f\in\cF_x$,
            \item[(ii)]   $\overline{\cU}\subseteq \bigcup_{f\in \cF_x} (\hat{D}_xf)^{-1}\cU$, \item[(iii)]  $m(D_xf)>1$, for $f\in\cF_x$.
        \end{itemize} 
        Then, there exists  $\rho>0$ such that for every $\tilde{\cF}\subseteq \diff_{\rm loc}^{1+}(M)$ 
         sufficiently close to $\cF$ in the $\cC^1$ topology,
 $\IFS(\tilde{\cF}\loc{V})$ is $\rho$-ergodic.
    \end{theorem} 
       
    Note that $V$ is not necessarily connected in Theorem \ref{thm:blender-Rd}. This yields a flexible tool to get a result similar to Theorem \ref{thm:blender} for every relatively compact open set in manifolds using local charts. 
    
    Let $M$ be a boundaryless manifold of dimension $d$ with a finite atlas of charts $\{(W_i,h_i)\}_{i\in I}$  on $M$ such that for any $i\in I$,  $h_i: W_i \to h_i(W_i)$ is a $\cC^2$ diffeomorphism  and $h_i(W_i)$ is an open disk in $\R^d$.
    Assume that 
    $h_i(W_i) \cap h_j(W_j)=\emptyset$ if $i\not=j$.
    For $i\in I$, let $W'_i$ be an open set  such that $\overline{W'_i}\subseteq W_i$ and  $M=\bigcup_{i\in I} W'_i$.
    Now, for an open set $V\subseteq M$ with compact closure, denote $V^* := \bigcup_{i\in I} h_i(V\cap W'_i)$.  Note that $V^*\subseteq\R^d$ is an open set with compact closure.
    Let $\cF \subseteq \diff_{\rm loc}^{1+}(M)$  and denote
    $\cF^*\subseteq \diff_{\rm loc}^{1+}(\R^d)$,
    \[\cF^* := \{h_j \circ f\circ h_i^{-1} {\gap}  f\in \cF\cup\{\id\}, ~i,j\in I ~\text{such that}~ x\in W_i ~\text{and}~ f(x)\in W_j \}.\]
    
    It is easy to see that the dynamical properties of $\cF$ are translated to the ones of $\cF^*$, and vice-versa.  Then, we get the following.

    \begin{theorem}\label{thm:blender-chart}
        Let $\cF, \cF^*, V,  V^*$ be as above. Let also $\cU \subseteq \sldr$. If $\cF^*, V^*, \cU$ satisfies the assumptions of Theorem \ref{thm:blender-Rd},  then there exists  $\rho>0$ such that for every $\tilde{\cF}\subseteq \diff_{\rm loc}^{1+}(M)$ 
        sufficiently close to $\cF$ in the $\cC^1$ topology, 
        $\IFS(\tilde{\cF}\loc{V})$ is $\rho$-ergodic. 
        
        In addition, if $M$ is compact,
        and the action of $\gen{\cF}$ is minimal, then it is $\cC^1$-stably ergodic and $\cC^1$-robustly minimal.
    \end{theorem}
    \begin{proof}
    The first part is an immediate consequence of Theorem \ref{thm:blender-Rd}, since the smooth maps send the sets of zero Lebesgue measure to the sets of zero Lebesgue measure. Moreover, Theorem \ref{thm:blender-Rd} shows that for every family $\tilde{\cF}$
     sufficiently close to $\cF$ in the $\cC^1$ topology,
    $\IFS(\tilde{\cF}\loc{V})$ is $\rho$-ergodic. 
    
    The second part is a duplication of the last parts of Theorem \ref{thm:blender-full}, in which minimality implies the stable covering of $M$ by the images of arbitrary balls of radius $\rho$ in $V$, under finitely many elements of $\gen{\cF}$.
    \end{proof}
    
    \section{Stably ergodic actions on manifolds}
    \label{sec:global}
    In this section, we prove Theorems \ref{thm:any-dim-sullivan}, \ref{thm:example} and \ref{thm:sphere-2} using the local results of the previous sections. 
    
    \subsection{Ergodic IFS on arbitrary manifold} In this subsection, we use the results of the previous section to construct a pair of diffeomorphisms generating a $\cC^1$-stably ergodic, $\cC^1$-robustly minimal IFS on any closed manifold $M$ of dimension $d$. Theorem \ref{thm:example} is a consequence of the following. 
    
    \begin{theorem}\label{thm:example-2}
    Given $s\in (1,+\infty]$, every closed Riemannian manifold $M$ admits a pair of diffeomorphisms in $\diff^{s}(M)$, generating a $\cC^1$-stably ergodic and $\cC^1$-robustly minimal semigroup action.
    \end{theorem}

    \begin{remark}\label{rem:no-ex}
    As far as the authors know, this is the first example of a stably ergodic action on a manifold of dimension greater than one.
    {In \cite{BFMS} and \cite{Sarizadeh} the existence of such actions on surfaces is claimed, however, with an argument that only works for conformal actions and mistakenly assumes ``conformality of maps with complex eigenvalues''. Indeed, it is easy to provide a non-quasiconformal sequence of matrices all with complex eigenvalues.}
    \end{remark}
     
    We will use the following simple lemma whose proof is very similar to the last part of the proof of Theorem \ref{thm:blender-full}. 
   
    \begin{lemma}\label{lem:local-to-global}
    Let $V$ be an open subset of compact manifold $M$ and $\cF$ be a family in $\diff^{1}(M)$. Assume that $\genplus{F}(V)=\langle{\mathcal{F}^{-1}}\rangle^+(V)=M$.
    If every measurable $\cF$-invariant  set has either zero or full Lebesgue measure in $V$, then $\ifs{F}$ is ergodic on $M$. Similarly, if for every $x\in M$, $\genplus{F}(x)$ is dense in $V$, then $\ifs{F}$ is minimal.   
    \end{lemma}
    \begin{proof}
    Let $S\subseteq M$ be measurable $\cF$-invariant set of positive Lebesgue measure, and let $x_0$ be a density point of $S$.  Since $\langle{\mathcal{F}^{-1}}\rangle^+(V)=M$,  one can find $n\in \N$ and $f_1,\ldots,f_n\in \cF$ such that $f_nf_{n-1}\cdots f_1(x_0)\in V$. This ensures that $S':=V\cap \langle \cF \rangle^+(S)$ has a positive Lebesgue measure in $V$. Since $\langle \cF\rangle^+(S')$ is $\cF$-invariant  and $\langle \cF\rangle^+(V)=M$, we deduce that $\langle \cF\rangle^+(S')$ equals $M$ up to a set of zero Lebesgue measure. Hence, $\ifs{\cF}$ is ergodic. 

    Proof of minimality is similar. Given $x_0\in M$, since  $\langle{\mathcal{F}^{-1}}\rangle^+(V)=M$,  one can find $n\in \N$ and $f_1,\ldots,f_n\in \cF$ such that $x':=f_nf_{n-1}\cdots f_1(x_0)\in V$. Now, since the closure of  $\langle \cF\rangle^+(x')$ contains $V$ and $\langle \cF\rangle^+(V)=M$, one deduce that the orbit of $x_0$ is dense in $M$ and so $\ifs{\cF}$ is minimal. 
    \end{proof}

    \begin{proof}[Proof of Theorem \ref{thm:example-2}] {The case $\dim(M)=1$ is known. It comes from Theorem \ref{thm:sullivan} and \cite{Gorodetski_Ilyashenko}. Indeed, Theorem \ref{thm:sphere} provides an explicit example for this case with independent proof. 
    
    So, we may assume that $\dim(M)\geq 2$. The same statement for robust minimality (without the ergodicity) is proved in \cite[Theorem A]{homburg_nassiri_2014}. Here, we carefully modify its proof and make use of our results in the previous sections to prove both robust minimality and stable ergodicity.  
    
    First, we establish some local construction on the Euclidean space and then in the next steps realize it by smooth diffeomorphisms on an arbitrary smooth manifold. Recall that we work with the $\cC^1$ topology on $\diff^s(M)$. }
    
    \medskip
   
    \noindent
    \textbf{Step 1.} \textit{Local construction: finitely many generators on $\R^d$.} 
    
    In this step, we construct a family of contracting affine transformations in $\R^d$ sufficiently close to Id satisfying the assumptions of Theorem \ref{thm:blender-cont}.

    Fix $\epsilon>0$ to be a sufficiently small number and let $\kappa \in (1,1+\epsilon)$.  
    Consider an open neighbourhood $\cU_0$ of the identity in $\sldr$ such that any element $D\in \cU_0$  is $\kappa$-conformal and by \eqref{eq:norm-4}, $\max\{\|D\|,\|D^{-1}\|\}<1+\epsilon$, hence $1-\epsilon<m(D)$. By Lemma \ref{lem:sldr-covering}, one can find a set $\cD\subseteq \cU_0$ consisting of $d^2$ elements, and an open set $\cU\subseteq \cU_0$ with 
    \begin{equation}\label{eq:cover-U}
    \overline{\cU}\subseteq \cD^{-1}\cU.
    \end{equation}
    For any $D\in \cD$ and $v\in \R^d$, define $T_{D,v}(x):=(1-\epsilon) D^{-1}(x)+v$.

     \begin{claim}\label{claim:V-J-delta}
         There exists $c(d)\in \N$, such that for any $\delta\in (0, \epsilon^2 ]$, there exists a finite set $J\subset V:=\Eucball{\delta}$ with $|J|=c(d)$ and 
    \begin{equation}\label{eq:cover-V}
    \overline{V}\subseteq\bigcup_{v\in J}T_{D,v}(V),
    \end{equation}
    for every $D\in \cD$.
     \end{claim}
    \begin{proof}
          First, note that there exists a finite set $J_0\subset \Eucball{1}$, such that $$\overline{\Eucball{1}}\subset \bigcup_{w\in J_0}(\Eucball{2/3}+w).$$ Let $c(d):=|J_0|$. Now, if $0\leq \delta \leq \epsilon^2$, then for $J:=\{\delta w:w\in J_0\}$,  
    \[\overline{V}\subset \bigcup\limits_{v\in J}(\Eucball{2\delta/3}+v).\] On the other hand, if $\epsilon$ is sufficiently small,  for every $D\in \cD$, $\Eucball{2\delta/3}\subset (1-\epsilon)D^{-1}V$.  Hence, for every $v\in \R^d$, $\Eucball{2\delta/3}+v\subset T_{D,v}(V)$. Therefore,
    \[\overline{V}\subset \bigcup\limits_{v\in J}(\Eucball{2\delta/3}+v) \subset \bigcup\limits_{v\in J} T_{D,v}(V).\]
    This finishes the proof of the claim.
    \end{proof}
    
Fix $\delta,J,V$ from Claim \ref{claim:V-J-delta}.
    We aim to show that the assumption of Theorem \ref{thm:blender-cont} are satisfied by taking $\cG:=\{T_{D,v}{\gap}D\in \cD, v\in J\}$, $\cW:=V\times \cU\subset \cE(\mathbb{R}^d)$, and $U:=\Eucball{1}\subseteq \R^d$.
    
   Note that every $T\in \cG$ is defined  everywhere in $\R^d$ and is  uniformly contracting, since $\|D_xT_{D,v}\|=(1-\epsilon)\|D^{-1}\|<1-\epsilon^2\leq 1-\delta$. From \eqref{eq:cover-V}, one gets that for any $D\in \cD$ and $x\in \overline{V}$, there exists $T\in \cG$ satisfying $y:=T^{-1}(x)\in V$ and $\pd_y T=D^{-1}$. This combined with (\ref{eq:cover-U}) implies that \[\overline{\cW}\subseteq \bigcup_{T\in \cG}\ef T(\cW).\] 
     So all the assumptions of Theorem \ref{thm:blender-cont} are satisfied for $\cG,\cW,U$.
Next, for $T\in \cG$, we define $\tilde{T}\in \diff^{\infty}(\R^d)$ sufficiently close to $T$ in the $\cC^1$ topology so that 
        \[\tilde{T}(x)=\begin{cases}
            T(x) & x\in \overline{\Eucball{1}}\\
            x & x\in \R^d\setminus \Eucball{2}
        \end{cases}\]
        Note that if $\epsilon$ is small enough, every  $T\in \cG$ is close to the identity in the $\cC^1$ topology. This implies all the elements of $\tilde{\cG}=\{\tilde{T}\;:\; T\in \cG\}$ are close to the identity in $\cC^1$ topology.  Now, it follows from  Theorem \ref{thm:blender-cont} that every $\tilde{\cG}\loc{U}$ invariant set of positive Lebesgue measure in $U$ contains $V$ up to a set of zero Lebesgue measure. Theorem \ref{thm:blender-cont} also guarantees that the orbit of every point $x\in U$ under $\tilde{\cG}\loc{U}$  is dense in $V$. Moreover, these properties are valid for every family that is sufficiently close to $\tilde{\cG}$ in the $\cC^1$ topology. 
        
    Enumerate the elements of $\tilde{\cG}$ by  $\tilde{T}_1,\ldots,\tilde{T}_{n_d}$, where $n_d:= d^2|J|=d^2c(d)$.

    \medskip
    \noindent
    \textbf{Step 2.} \emph{Local ergodicity: a pair of generators on M.}

    First note that for every closed manifold admits a
    $\cC^\infty$ Morse-Smale diffeomorphism $f$ with a unique attracting periodic orbit $\orbit_f(p)$ with arbitrarily large period $N$  so that $f^N$ is close to the identity. Moreover, we may assume that all the periodic points of $f$ are fixed except those in $\orbit_f (p)$. Such a diffeomorphism $f$ is explicitly constructed in \cite{homburg_nassiri_2014} by deforming the time one map of the gradient flow of a suitable Morse function with a unique minimum point.

   We choose such $f$ with $N>2n_d$. Denote the set of all other periodic points of $f$ by $Q_f$, which is a finite set and is equal to $\mathrm{Fix}(f)$.  Let  $U_0\subseteq M$ be a small neighbourhood of $p$ satisfying
    \begin{itemize}
     \item $U_0$ is contained in the domain of attraction of $p$. In particular, the orbit of $U_0$ under $f$ does not intersect $Q_f$.
     \item for $i=0,\ldots, N-1$, the sets $U_i$ are pairwise disjoint and  $U_N\subseteq U_0$, where $U_i:=f^i(U_0)$.    
    \end{itemize}    

     Notice that from Step 1,  $V=\Eucball{\delta}\subset U=\Eucball{1}$ are open sets in $\R^d$, while  $U_0,\ldots,U_N$ are open subsets of $M$.

    Now, choose a diffeomorphism $\phi:U_0\to \R^d$ with $\phi(p)=0$ and $\overline{\Eucball{2}}\subseteq {\phi(U_N)}$. Then, define $h\in \diff^\infty(M)$ as follows (see Figure \ref{fig:thm-example}).
    
    \begin{itemize}
        \item[-] $h$ equals to the identity on  $M\setminus\bigcup_{i=1}^{2n_d}U_i$.
        \item[-] For $1\leq i \leq n_d$,
\begin{equation}
    \label{eq:definition_of_h|U_i}
    \begin{aligned}
 h|_{U_{2i-1}}& :=f^{2i-1}\circ \phi^{-1}\circ \tilde{T}_i\circ \phi \circ f^{-2i+1},\\
         h|_{U_{2i}}& :=f^{-1}\circ (h|_{U_{2i-1}})^{-1} \circ f,
    \end{aligned}
\end{equation}
        
        where 
$\tilde{T}_1,\ldots,\tilde{T}_{n_d}$ are the elements of $\tilde{\cG}$ defined in Step 1.
        \end{itemize}

It is clear that for every $i=1,\ldots,n_d$, $(h\circ f^{-1})^2|_{U_{2i-1}}=f^{-2}|_{U_{2i-1}}$ This implies the following equality which will be used in Step 3. 
\begin{equation}\label{eq:hf^{-1}N}
     (h\circ f^{-1})^{-N}|_{U_0}=f^{N}|_{U_0}.
\end{equation}

    \begin{claim}\label{claim-local-erg}
    If $f$ is sufficiently close to the identity, then for every    $\tilde{f},\tilde{h}\in \diff^s(M)$  sufficiently $\cC^1$-close to $f,h$, respectively, we  have the following properties: 
    \begin{itemize}
        \item[(i)]  Every measurable $\langle \tilde{f},\tilde{h}\rangle^+$-invariant set with positive Lebesgue measure in  $V_0:=\phi^{-1}(V)\subseteq M$ contains $V_0$ up to a set of zero Lebesgue measure.  
        \item[(ii)]  For every $x\in V_0$, $\langle \tilde{f},\tilde{h}\rangle^+(x)$ is dense in $V_0$. 
    \end{itemize}
       
    \end{claim}
    \begin{proof}
  The idea is to use the construction of $h$ to reduce the problem to that of Step 1. 

  Define $h_i:U_0\to U_0$ by $h_i:=\phi^{-1}\circ \tilde{T}_i\circ \phi$. Hence, the elements of the family $\cH:=\{h_1,\ldots,h_{n_d}\}$ are simultaneously smoothly conjugated to the elements of $\tilde{\cG}=\{\tilde{T}_1,\ldots,\tilde{T}_{n_d}\}$.
    So, it is enough to show that if $\tilde{f},\tilde{h}$ are sufficiently $\cC^1$-close to $f,h$ respectively, then  $\langle \tilde{f},\tilde{h}\rangle^+$ contains a family $\tilde{\cH}$   sufficiently $\cC^1$-close to $\cH$. This would finish the proof of Claim \ref{claim-local-erg} by the conclusion of Step 1. 
    
    For every $1\leq i \leq n_d$, $h_i\circ f^N$ maps $U_0$ into $U_0$, and if $f^N$ is sufficiently close to the identity, $f^N\circ h_i|_{U_0}$ is sufficiently close to $h_i$ in the $\cC^1$ topology.
    On the other hand, it follows from \eqref{eq:definition_of_h|U_i} that for every $1\leq i \leq n_d$ and every $x\in U_0$,
     \begin{equation}\label{eq:h_i_fN}
         f^{N-2i+1}\circ h\circ f^{2i-1}(x)=f^N\circ h_i(x),
     \end{equation}
    Now, for every  small $\cC^1$ perturbation $\tilde{f},\tilde{h}$ of $f,h$,  respectively, define \[\tilde{\cH}:=\{\tilde{f}^{N-2i+1}\circ \tilde{h}\circ \tilde{f}^{2i-1}|_{U_0}\;:\; 1\leq i\leq n_d\}.\]
    It is clear that $\tilde{\cH}\subseteq \langle\tilde{f},\tilde{h}\rangle^+$ and by \eqref{eq:h_i_fN} it is $\cC^1$-close to $\cH$.    
    \end{proof}

    \begin{remark}\label{rm:large-perturbation}
       The arguments of Step 2 rely on the definition of $f,h$ inside $\bigcup_{i=0}^NU_i$. Thus,  the conclusion of Claim  \ref{claim-local-erg} remains true
        for any $\hat{f},\hat{h}\in \diff^1(M)$ that equals to $f,h$ on $\bigcup_{i=0}^NU_i$, respectively. 
    \end{remark}

    \begin{figure}[t]
        \centering
    \includegraphics[width=.55\textwidth]{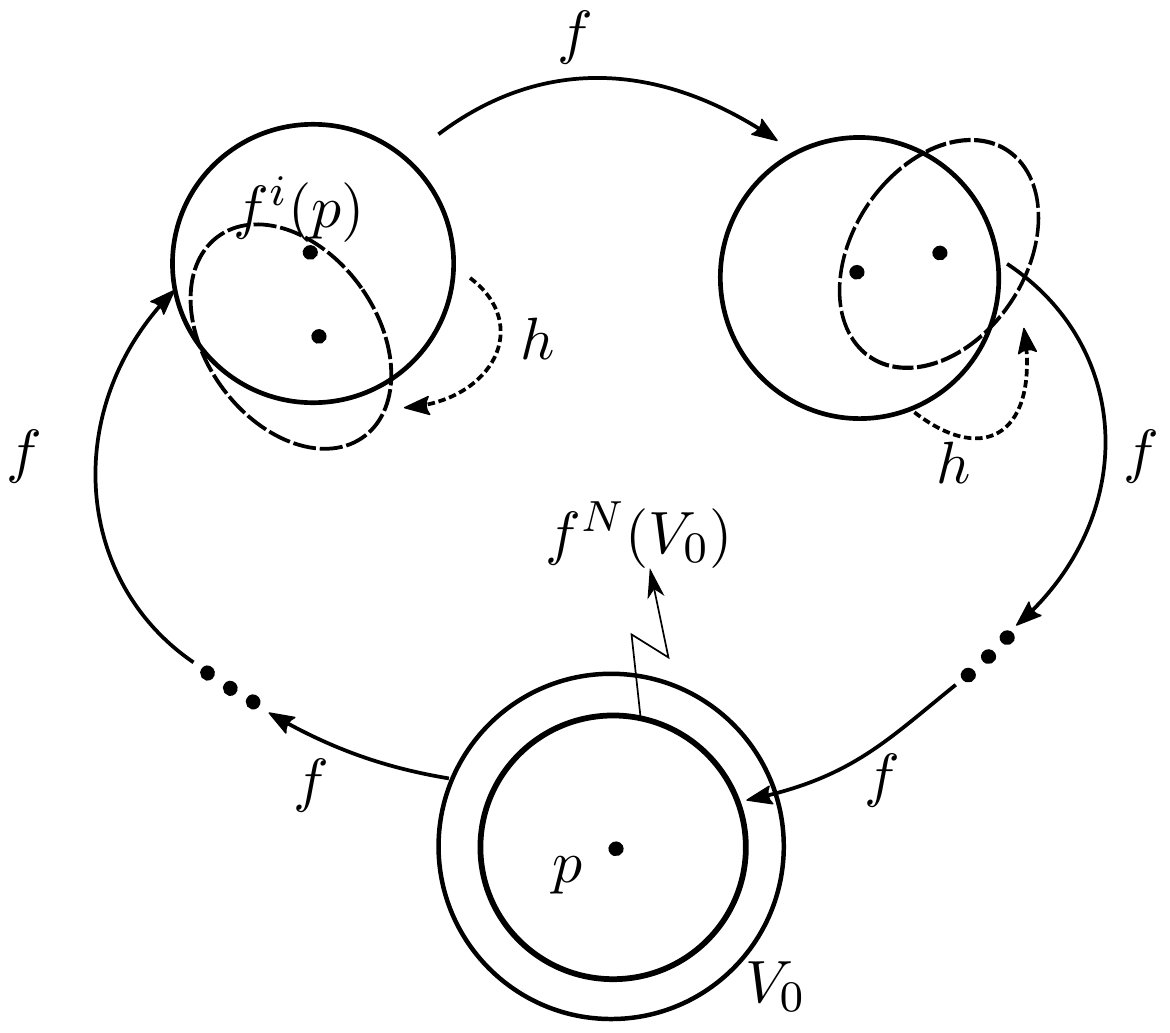}
        \caption{Local constructions of $f$ and $h$.}
        \label{fig:thm-example}
    \end{figure}
    
    \medskip
    \noindent
    \textbf{Step 3.} \emph{Global construction: a pair of generators on M.}

   Take diffeomorphisms $f,h$ given in Step 2.
    Recall that $f$ is a Morse-Smale diffeomorphism with a unique attracting periodic orbit $\orbit_f(p)$. The forward orbit of every point under $f$ converges either to $\orbit_f(p)$ or to an element of $Q_f$. Therefore, 
    \begin{equation}\label{eq:f-decomposition-M}
        M=\bigcup_{q\in Q_f}W_f^s(q)\cup W_f^s(\orbit_f(p)),
    \end{equation}
    where $W_f^s$ denotes the stable manifold of $f$. Since $\orbit_f(p)$ is the unique attracting periodic orbit for $f$, $W_f^s(\orbit_f(p))=\bigcup_{i>0} f^{-i}(V_0)$ is an open and dense subset of $M$, where $V_0$ is defined in Claim \ref{claim-local-erg}.

Now, pick $\psi\in \diff^\infty(M)$ 
such that
    \begin{itemize}
        \item[- ] 
        $\psi(Q_f)\cup \psi^{-1}(Q_f)\subseteq 
        W^s_f(\orbit_f(p))$. 
        \item[- ] $\psi$ equals the identity outside a sufficiently small neighbourhood  $W$ of $Q_f$ with $W\subset M\setminus \cup_{i=0}^N \overline{U_i}$. This, in particular, implies  $h\circ \psi=\psi\circ h$. 
    \end{itemize}

 Finally, let $g:=\psi \circ (h \circ f^{-1})\circ \psi^{-1}\in \diff^{\infty}(M)$ and $\cF:=\{g,f\}$.

Notice that,
$g= \psi \circ f^{-1}\circ \psi^{-1}$ on the set $M\setminus\bigcup_{i=0}^{2n_d}U_i$.
Thus, the points in $Q_g:=\psi(Q_f)$ are hyperbolic fixed points of $g$.
On the other hand, by \eqref{eq:hf^{-1}N}, $p$ is the unique periodic point of $g$ in $U_0$, and it is repelling. This means that the nonwandering set of $g$ consists of hyperbolic periodic points and it is equal to $Q_g\cup \orbit_g(p)$. So,
\begin{equation}\label{eq:g-decomposition-M}
        M=\bigcup_{q\in Q_g}W_g^u(q)\cup W_g^u(\orbit_g(p)),
    \end{equation}
where, $W_g^u$ denotes the unstable manifold of $g$. Moreover, $\psi^{-1}(Q_f)\subseteq W^s_f(\orbit_f(p))$ implies that 
\begin{equation}\label{eq:Q_f}
    Q_f\subseteq W^u_g(\orbit_g(p)).
\end{equation}

We claim that $\IFS(\cF)$ is $\cC^1$-stably ergodic and $\cC^1$-robustly minimal on $M$. 

First, observe that by the definition of $g$, $g\circ f=h$ on $\cup_{i=0}^{N}U_i$. Hence, according to Remark \ref{rm:large-perturbation},  we deduce that the local ergodicity and minimality obtained in Step 2 hold for $\langle g\circ f , f\rangle^+$ and also for $\langle f,g\rangle^+$ as  $\langle g\circ f , f\rangle^+\subseteq \langle f,g\rangle^+$. Now, in view of  Lemma \ref{lem:local-to-global}, it suffices to prove the following.

 \begin{claim} 
        For every $\tilde{f},\tilde{g}$ sufficiently $\cC^1$-close to $f,g$, respectively. 
        \[\langle \tilde{f},\tilde{g}\rangle^+(V_0)=\langle \tilde{f}^{-1},\tilde{g}^{-1}\rangle^+(V_0)=M.\]
    \end{claim}
    \begin{proof}
        We first prove the statement for $f,g$, themselves. If $x\in W^s_f(\orbit_f(p))$, then there exists $m\in \N$ with $f^m(x)\in V_0$. Otherwise, if $x\in \bigcup_{q\in Q_f}W^s_f(q)$, then there exists $q\in Q_f$ such that $f^n(x)\to q$, as $n$ tends to infinity.
        On the other hand,  $g^n(q)\to q'$, as $n$ tends to infinity, for some  $q'\in Q_g=\psi(Q_f)\subseteq W^s(\orbit_f(p))$. Therefore, for sufficiently large $n_1,n_2,n_3\in \N$, we have  
        $f^{n_3}\circ g^{n_2}\circ f^{n_1}(x)\in V_0$. 
        In view of \eqref{eq:f-decomposition-M}, this implies that 
        $\langle f^{-1},g^{-1}\rangle^+(V_0)=M$. 
 Analogously, one can deduce from 
 \eqref{eq:g-decomposition-M} and \eqref{eq:Q_f}   that $\langle {f},{g}\rangle^+(V_0)=M.$

         Finally,  Since $V_0$ is open and $M$ is compact, there is a finite subset $\cF_0\subseteq \langle f,g\rangle^+ $ such that $\cF_0(V_0)=\cF_0^{-1}(V_0)=M$. Since small perturbations of $f,g$ imply small perturbations of elements of $\cF_0$, we deduce that the conclusion of the claim holds for every  $\tilde{f},\tilde{g}$ sufficiently $\cC^1$-close to $f,g$. 
  
    \end{proof}

    Hence, the proof of Theorem \ref{thm:example-2} is complete. 
    \end{proof}
    
    \subsection{Proof of Theorem \ref{thm:any-dim-sullivan}} 
    We present a variant of Theorem \ref{thm:any-dim-sullivan} for diffeomorphisms between open sets of $M$, i.e., they are not necessarily globally defined on the manifold.
    In particular, it implies Theorems \ref{thm:any-dim-sullivan}.
    
    \begin{theorem}\label{thm:local-sullivan}
    Let $M$ be a closed smooth Riemannian manifold and 
       $\cF\subseteq \diff_{\rm loc}^{1+}(M)$. Assume that for any $(x,v)\in T^1M$ there exists $f\in  \cF$ such that $m(D_xf) > 1$ and  $\|\hat{D}_xf|_{v^{\perp}} \|<1$. Then, there exists $\rho>0$ such that the action of  $\gen{\cF}$ is stably $\rho$-ergodic.
       
       Moreover, if the action of $\gen{\cF}$ is minimal, it is $\cC^1$-stably ergodic in $\diffloc{1+}(M)$ and is $\cC^1$-robustly minimal.
    \end{theorem}
    \begin{proof}
    The proof follows from Lemma \ref{lem:direc-cont} and  Theorem \ref{thm:blender-full}. By the compactness of $T^1M$, one can find a finite subset $\cF_0\subseteq \cF$
    such that for any $(x,v)\in T^1M$  there exists $f\in \cF_0$ with $m(D_xf)>1$ and $\|\pd_x f|_{v^\perp}\|<1$. For any $f\in\cF_0$, denote $U_f:=\{x\in M{\gap}m(D_xf)>1\}$. Now, by Lemma \ref{lem:direc-cont}, the assumptions of Theorem \ref{thm:blender-full} are satisfied for family $\hat{\cF}_0:=\{f|_{U_f}{\gap}f\in \cF_0\}$. Therefore, there exists $\rho>0$ such that $\IFS(\hat{\cF}_0)$ is stably $\rho$-ergodic. In particular,  $\IFS(\cF)$ is stably $\rho$-ergodic. 
   
    Then, similar to the proof of the last part of Theorem \ref{thm:blender-full}, we use minimality of the action of $\gen{\cF}$ to stably cover $M$ by the images of an arbitrary ball of radius $\rho$ under a finite set $\cF_1\subseteq \gen{\cF}$. This proves stable ergodicity and robust minimality of the action of $\gen{\cF}$.
    \end{proof}

\subsection{Proof of Theorem \ref{thm:sphere-2}}
     We first prove the following special case of Theorem \ref{thm:sphere-2}.
     
    \begin{theorem}\label{thm:sphere}
    Let $d\geq 1$ and $\{A_1,\ldots,A_k\}\subseteq \mathrm{SO}(d+1)$ generates a dense subgroup of $\mathrm{SO}(d+1)$. Then, for any  $A_0 \in \mathrm{SL}(d+1,\R)\setminus\mathrm{SO}(d+1)$, the {natural} action of the group generated by $\{A_0,\ldots,A_k\}$ on $\SS^d$ is $\cC^1$-stably ergodic in $\diff^{1+\alpha}(\SS^d)$. Moreover, it is $\cC^1$-robustly minimal.
    \end{theorem} 
    
    {We denote by $f_A$, the action of $A\in \mathrm{SL}(d+1,\R)$ on $\sd$, which is defined by $x\mapsto \frac{Ax}{|Ax|}$. Also, denote the standard orthonormal basis of $\R^{d+1}$ by $(e_1,\ldots,e_{d+1})$. One can  easily check that if  the diagonal matrix $\hat{A}=\diag(r_1,\ldots,r_{d+1}),$ satisfies 
    \[0<r_{d+1}<\min\limits_{1\leq i \leq d}r_i,\text{and}~ \max\limits_{1< i \leq d} r_i^d<r_1\cdots r_d,\]
    then
    \begin{equation}\label{eq:direc-cont-sphere}
        m(D_{e_{n+1}}f_{\hat{A}})>1, \text{and}~ \|\pd_{e_{d+1}}f_{\hat{A}}|_{W}\|<1,
    \end{equation}
    where $W$ is the $(d-1)$-dimensional subspace in $T_{e_{d+1}}\sd$ perpendicular to $e_1$. Note that the inequalities in (\ref{eq:direc-cont-sphere}) are stable under perturbations of $e_{d+1}$, $W$ and $f$. More precisely,
    \begin{itemize}
       \item[($*$)] There exist $\lambda<1$, $\epsilon>0$ and open neighbourhoods $U$ of $e_{d+1}$ such that for any $\cC^1$ map $f$ sufficiently close to $f_{\hat{A}}$ in the $\cC^1$ topology, and any  $(p,v)\in T^1\sd$ with $p\in U$ and $\angle(v,e_1)<\epsilon$, one has $\|\pd_p f|_{v^\perp}\|<\lambda$ and $\|D_pf\|>\lambda^{-1}$.
    \end{itemize}}
    We will need the following {linear algebraic} lemma. 
    \begin{lemma}\label{lem:linear-algebra}
    Given any matrix $D\in \mathrm{SL}(d+1,\mathbb{\R})\setminus \SOd$, there exist $n>0$, $R_0,R_1,\ldots,R_n\in \SOd$, and $\alpha_1,\ldots,\alpha_n\in\{-1,+1\}$ such that $$R_0D^{\alpha_1}R_1D^{\alpha_2}R_2\ldots R_{n-1} D^{\alpha_n} R_n=\diag(r_1,\ldots,r_{d+1}),$$ where  $0<r_{d+1}<\min\limits_{1\leq i \leq d}r_i$ and $\max\limits_{1< i \leq d} r_i^d<r_1\cdots r_d$.
    \end{lemma}
    \begin{proof}
    For any permutation $\sigma$ of $\{1,\ldots,d+1\}$, there are signs  $\varepsilon_1,\ldots,\varepsilon_{d+1}\in\{-1,+1\}$ such that the linear map $R$  defined by  \[R(x_1,\ldots,x_{d+1})=(\varepsilon_1 x_{\sigma(1)},\ldots,\varepsilon_{d+1} x_{\sigma(d+1)}),\] 
    is an element of $\SOd$. Denote an element $R$ associated to $\sigma$ in this way by $R_\sigma$. 
  
    By means of the singular value decomposition, one can write $D=RD'R'$ with $R,R'\in \SOd$ and $D'$ diagonal. So, it suffices to prove the lemma for diagonal matrices $D$. By replacing $D$ with some $R_\sigma D R_{\sigma}^{-1}$, if necessary, suppose that $D=\diag(s_1,\ldots,s_{d+1})$ with $0< |s_{d+1}|\leq |s_d|\leq \cdots \leq  |s_1|$. Since $D\not\in \SOd$, $|s_1|>1$. Let $\Sigma$ be the set of all permutations $\sigma$ of $\{1,\ldots,d+1\}$ with $\sigma(1)=1$. 
    Then, $D_1:=\prod_{\sigma\in \Sigma} R_\sigma DR_\sigma^{-1}=\diag(t_1,\ldots,t_{d+1})$ satisfies
    \[|t_1|>|t_2|=\cdots=|t_{d+1}|=1.\]
    Let $\sigma_0$ be the permutation on $\{1,\ldots,d+1\}$ with $\sigma_0(i)=d+1-i$. Hence,
    \[R_{\sigma_0}D_1^{-1}R_{\sigma_0}^{-1}=\diag(t_{d+1}^{-1},\ldots, t_1^{-1}),\]
    and so $D_2=D_1R_{\sigma_0}D_1^{-1}R_{\sigma_0}^{-1}=\diag(r_1,\ldots,r_{d+1})$ with $|r_1|>1>|r_{d+1}|$ and $|r_2|=\cdots=|r_d|=1$. Now, $D_2^2$ has positive diagonal entries and satisfies the conditions. 
    \end{proof}
   
    \begin{proof}[Proof of Theorem \ref{thm:sphere}]
    Denote $\cF:=\{f_{A_1},\ldots,f_{A_k}\}$ and $\cF_0:=\cF\cup \{f_{A_0}\}$. Note that ergodicity and minimality of the action of $\gen{\cF}$ is {guaranteed, since} $\overline{\gen{A_1,\ldots,A_k}}=\SOd$. So, we only need to show the stability under $\cC^1$ perturbations. To this end, we use the minimality of the isometric action of $\gen{\cF}$ to transfer the properties of the derivatives in some open set to the whole manifold. 
  
    \medskip
    \noindent
    \textbf{Claim} (Reduction to Theorem \ref{thm:local-sullivan}){\bf.}\emph{ There exists a finite set $\cF_1\subseteq \langle\cF_0\rangle$ such that for any $(x,v)\in T^1\sd$, $m(D_xg)>1$ and $\|\pd_xg|_{v^\perp}\|<1$ for some  $g\in\cF_1$.}
    \begin{proof}
    The fact $\overline{\gen{A_1,\ldots,A_k}}=\SOd$ combined with Lemma \ref{lem:linear-algebra} for $D=A_0$ implies that there exists $A\in \gen{A_0,A_1,\ldots,A_k}$ sufficiently close to $\hat{A}$, defined above, such that $(*)$ holds for $f_A$ and its $\cC^1$ perturbations. Fix this $A$ for the rest of the proof and denote $f:=f_A$, for simplicity.
   
    {It follows from the minimality of the action on $T^1\sd$ that for any $(x,v)$, there exists $h\in \gen{\cF}$ such that for every $(y,w)$ in a neighbourhood of $(x,v)$, $h(y)\in U$ and  $\angle(D_yh(w),e_1)<\epsilon$.  Now,  by $(*)$, $g:=f\circ h\in \gen{\cF_0}$ satisfies $m(D_yg)>\lambda^{-1}$ and $\|\pd_yg|_{w^\perp}\|<\lambda$. Finally, by the compactness of $T^1\sd$, one can choose a finite subset $\cF_1$ of $\gen{\cF_0}$ satisfying these properties.} 
    \end{proof}
    Since the action of $\cF_0$ on $\sd$ is minimal, this claim combined with Theorem \ref{thm:local-sullivan} implies that the action of $\gen{\cF_0}$ is $\cC^1$-stably ergodic and $\cC^1$-robustly minimal in $\diff^{1+}(\SS^d)$. 
    \end{proof} 

    \begin{proof}[Proof of Theorem \ref{thm:sphere-2}]
    Consider a finite family $\{A_1,\ldots,A_k\}\subseteq \SOd$ generating a dense subgroup of $\SOd$. The existence of such elements for $d=1$ is trivial and for $d\geq 2$,  is granted by  \cite{Kuranishi}, as $\SOd$ is a semi-simple Lie group. Since $\gen{\cF}$ is dense in $\SOd$, it contains elements $\tilde{A}_1,\ldots,\tilde{A}_k$ arbitrarily close to $A_1,\ldots,A_k$. On the other hand, $\cF$ has an element $A_0\in \mathrm{SL}(d+1,\R)\setminus\SOd$. So, Theorem \ref{thm:sphere} implies that the natural action of $\gen{\cF}$ on $\sd$ is stably ergodic and robustly minimal. \end{proof}

\begin{proof}[Proof of Corollary 
\ref{cor:sphere-2}]
    As mentioned before,  it is known that for $d\geq 2$, every generic pair $A, B$ of elements near the identity generates a dense subgroup of  $\sldr$ \cite{Kuranishi}. So,  $\gen{A,B}$ satisfies the assumptions of Theorem \ref{thm:sphere-2}. This  completes the proof.
\end{proof}

\begin{example}\label{ex:non-ergodic}
    Let $A, B$ be elements of $\mathrm{SL}(2,\R)$ in some neighbourhoods of $(\begin{smallmatrix} 100 & 0  \\ 0 & 0.01 \end{smallmatrix})$ and $(\begin{smallmatrix} 10 & 11  \\ 9 & 10 \end{smallmatrix})$, respectively. By the ping-pong argument one can easily show that the action of $\gen{A,B}$ on $\SS^1$ is neither minimal nor ergodic.
\end{example}

    \subsection{Proof of Corollary \ref{cor:stationary}}\label{sec:proof-stationary}

    Recall that given a finitely supported probability measure $\nu$ on  $\diff^{1+\alpha}(M)$, a probability measure $\mu$ on $M$ is called {\it $\nu$-stationary},  if 
    \(\sum \nu(\{f\}) f_*\mu =\mu\), where the sum is over the support of $\nu$
    and $f_*\mu$ is the pushforward of the measure $\mu$ under $f$. The set of $\nu$-stationary probability measures form a compact convex subset of probability measures on $M$. The {\it ergodic} stationary measures are extreme points of this convex set.

\begin{proof}[Proof of Corollary 
\ref{cor:stationary}]
 Let $\cF$ be the support of the measure $\nu$  and for every  $f\in \cF$ denote $\nu_f:=\nu(\{f\})>0$. Suppose that $\mu$ is an ergodic $\nu$-stationary measure on $M$. We aim to prove that if the measures $\mu$ and $\leb$ are not singular, then they are equivalent. So we assume that $\mu$ is not singular with respect to $\leb$.  
 
 First, we show that $\leb$ is absolutely continuous with respect to $\mu$. Assume to the contrary that there exists a measurable set $S\subseteq M$ with $\mu(S)=0$ but $\leb(S)>0$. Since $\mu$ is $\nu$-stationary,
\[0=\mu(S)=\sum\limits_{f\in \cF} \nu_f\cdot f_*\mu(S)=\sum\limits_{f\in \cF} \nu_f \cdot\mu(f^{-1}(S))\]
Therefore, we deduce that for every  $f\in \cF$,  $\mu(f^{-1}(S))=0$. Similarly, one can prove that for every $f\in \langle\cF^{-1}\rangle^+$, $\mu(f^{-1}(S))=0$. Therefore, $\mu(S^-)=0$, where $S^-:=\bigcup_{f\in \langle \cF^{-1}\rangle^+}(S)$. On the other hand,  $\leb(S^-)>0$ and since $M\setminus S^-$ is an $\cF$-invariant measurable set, by the ergodicity of the action of semigroup $\langle\cF\rangle^+$ on $M$, we have $\leb(M\setminus S^-)=0$. Therefore, $\leb(S^-)>0$ which by absolute continuity of $\leb$ with respect to $\mu$ implies that $\mu(S^-)>0$, which is a contradiction.

Next, we only need to show that $\mu$ is absolutely continuous with respect to $\leb$. By Lebesgue's decomposition theorem, there exists probability measures $\mu_{a}$ and $\mu_{s}$ which are absolutely continuous and singular with respect to $\leb$, respectively, and $t\in [0,1]$ such that $\mu=t\mu_{a}+(1-t)\mu_{s}$. 
   \begin{claim*}   
 $\mu_{a}$ and $\mu_{s}$ are both $\nu$-stationary.
    \end{claim*} 
    \begin{proof}
    Denoting  $\tilde{\mu}_{a}:=\sum_{f\in \cF} \nu_f\cdot f_*\mu_{a}$ and $\tilde{\mu}_{s}:=\sum_{f\in \cF} \nu_f\cdot f_*\mu_{s}$, we need to show that $\tilde{\mu}_{a}=\mu_{a}$ and $\tilde{\mu}_{s}=\mu_{s}$. First, we show that $\tilde{\mu}_{a}$ is absolutely continuous with respect to $\leb$. Indeed, if $\leb(S)=0$, then for every $f\in \cF$, $\leb(f^{-1}(S))=0$ and by absolute continuity of $\mu_{a}$ 
    \[\sum\limits_{f\in \cF}\nu_f \cdot \mu_{a}(f^{-1}(S))=0.\]
    On the other hand, we will show that $\tilde{\mu}_{s}$ is singular to $\leb$. Since $\mu_{s}$ is singular to $ \leb$, there exists a measurable set $A$ with $\leb(A)=\mu_{s}(A^c)=0$. Clearly, for every diffeomorphism $h$,  $\leb(h(A))=0$ and so $\leb(B)=0$, where 
    \[B:=\bigcup\limits_{h\in \cF} h(A)\]
    
    This implies  
    \begin{align*}
        \tilde{\mu}_{s}(B^c)& =\sum\limits_{f\in \cF} \nu_f\cdot  f_*\mu_{s}(B^c) = \sum\limits_{f\in \cF} \nu_f\cdot  f_*\mu_{s} \Big(\bigcap_{h\in \cF}h(A^c)\Big)\\
        & \leq \sum\limits_{f\in \cF} \nu_f\cdot  f_*\mu_{s}\big(f(A^c)\big)
        =\sum\limits_{f\in \cF} \nu_f\cdot  \mu_{s}(A^c)=0.
    \end{align*}
    Therefore, $\tilde{\mu}_{s}(B^c)=\leb(B)=0$ implying  that $\tilde{\mu}_{s}$ is singular with respect to $\leb$. 
    Now, since $\mu=t\mu_{a}+(1-t)\mu_{s}$, one has  
    \[\mu=\sum\limits_{f\in \cF} \nu_f\cdot  f_*\mu=t\Big(\sum\limits_{f\in \cF} \nu_f\cdot  f_*\mu_{a}\Big)+(1-t) \Big(\sum\limits_{f\in \cF} \nu_f\cdot  f_*\mu_{s}\Big)=t\tilde{\mu}_{a}+(1-t)\tilde{\mu}_{s}.\]
    Due to the uniqueness of the decomposition provided by Lebesgue's decomposition theorem, we get  $\tilde{\mu}_{a}=\mu_{a}$ and $\tilde{\mu}_{s}=\mu_{s}$, that is, both measures $\mu_{a}$ and $\mu_{s}$ are $\nu$-stationary.   \end{proof}
Finally, since $\mu$ is an extreme point in the space of stationary probability measure, one gets either $\mu=\mu_{a}$ or $\mu=\mu_{s}$. The latter case is impossible as we have assumed that $\mu$ is not singular to $\leb$. Therefore, $\mu=\mu_{a}$ is absolutely continuous with respect to $\leb$.   
    This finishes the proof of Corollary \ref{cor:stationary}.
\end{proof}

    \section{Some questions}\label{sec:questions}
 
    \subsection{The number of generators for stably ergodic actions}  \label{subsec:opt} 
 Theorem \ref{thm:example} states that every manifold $M$ admits a stably ergodic semigroup action generated by two diffeomorphisms. A natural question to ask is whether 
 the number of generators in Theorem \ref{thm:example} is optimal. In other words, 
    \begin{question}
    Does there exist a manifold $M$ with a stably ergodic (w.r.t.  Leb.) diffeomomorphism in $\diff^{1+\alpha}(M)$?
    \end{question}
    Two observations support a negative answer to this question. First, {\it no Anosov diffeomorphism is $\cC^1$-stably ergodic (w.r.t.  Leb.) in $\diff^s(M)$ ($s\geq 1$)}. Second, some manifolds do not admit  stably transitive diffeomorphisms. More precisely, {\it there is no $\cC^1$-stably ergodic cyclic group in $\diff^{1+\alpha}(M)$, if $M$ is either the circle, a closed surface, a 3-manifold that does not admit partially hyperbolic diffeomorphisms (e.g. $\SS^3$), or a manifold whose tangent bundle does not split (e.g. $\SS^{2k}$).}

    It follows from \cite{Gurevich-Oseledets} that any $\cC^2$ Anosov diffeomorphism, which is ergodic with respect to Lebesgue, admits a unique invariant measure in the Lebesgue class. 
    On the other hand, by \cite{Livsic-Sinai}, the set of $\cC^2$ Anosov diffeomorphisms admitting no absolutely continuous invariant measure form an open and dense subset in the $\cC^1$ topology (satisfying an explicit condition on the derivative of some periodic point). Since every Anosov diffeomorphism can be approximated by $\cC^2$ Anosov diffeomorphisms, it follows that no Anosov diffeomorphism is stably ergodic (w.r.t.  Leb.) in $\diff^{1+\alpha}(M)$.

    Next, it is a consequence of  \cite{Mane}, \cite{Diaz_Pujals_Ures} and \cite{Bonatti_Diaz_Pujals} that some forms of hyperbolicity (and thus splitting of the tangent bundle) can be obtained from robust transitivity in  $\diff^{1}(M)$. Moreover, on $\T^2$,  a $\cC^1$-robustly transitive diffeomorphism is indeed an Anosov diffeomorphism.  As a matter of fact, the same proofs work if one considers the $\cC^1$ topology in the space $\diff^{1+\alpha}(M)$, as we do here. Thus, none of the manifolds listed above (except $\T^2$) do admit a stably transitive diffeomorphism. 
    The claims on $\SS^3$ and $\SS^{2k}$ follow from \cite{Brin-Burago-Ivanov_2004} and the obstruction theory in topology \cite{Milnor-Stasheff}, respectively. 
 
    These arguments raise the following questions (see also \cite{Avila_Bochi}).

    \begin{question}
    Which ergodic partially hyperbolic diffeomorphisms in $\diff^2(M)$ admit an invariant probability measure equivalent to the Lebesgue measure?
    \end{question}
     
    \begin{question}
    Is it true that a generic diffeomorphism in $\diff^2(M)$ admits no invariant measure in the Lebesgue class?
    \end{question} 
     
    \subsection{Ergodicity vs. quasi-conformality} 
    In the theory of stable ergodicity in $\diff^{2}_{\rm vol}(M)$, usually, the ergodicity follows from (some forms of) hyperbolicity using a Hopf type argument. Clearly, any form of hyperbolicity obstructs  the existence of quasi-conformal orbits. In contrast, as discussed above, even Anosov diffeomorphisms are not stably ergodic (w.r.t.  Leb.) in $\diff^{1+\alpha}(M)$. Moreover, the existence of quasi-conformal orbits is a crucial ingredient of our arguments to establish stable (local) ergodicity. One may ask the following.  
    
    \begin{question}\label{q:q.c-branch}
    Does there exist a $\cC^1$-stably ergodic finitely generated semigroup in $\diff^2(M)$ such that (stably) the set of points having quasi-conformal orbit-branches has zero Lebesgue measure?
    \end{question}

    In Theorem \ref{thm:any-dim-sullivan}, while the contraction hypothesis for the normalized derivative on hyper-subspaces obstructs conformality, it allows one to obtain (stably) quasi-conformal orbit-branches at every point. Inspired by the work of \cite{Avila_Bochi_Yoccoz} on semigroups of ${\rm SL}(2, \R)$ we ask the next question concerning the optimality of this assumption in dimension 2.

    \begin{question}\label{q:covering}
    Let $M$ be a closed surface, $\cF\subseteq \diff^1(M)$ be  finite. Suppose that the action of $\genplus{F}$ has (stably) quasi-conformal orbit-branches at every point. Is it true that for $s\geq 1$, $\cF$ can be $\cC^s$-approximated by the families $\tilde{\cF}$ satisfying the following condition?
    \begin{itemize}
        \item For every $(x,v)\in T^1M$, there exists $f\in \langle \tilde{\cF}\rangle^+$ with $|\pd_xf(v)|<1$.
    \end{itemize}
    \end{question}
       
    \subsection{Stationary measures}\label{subsec:st-meas}

    It is clear that the existence of an ergodic stationary measure in the measure class of $\leb$ for some distribution on a semigroup $\genplus{F}$  implies ergodicity of $\IFS(\cF)$ (w.r.t. $\leb$). Conversely, ergodicity with respect to the Lebesgue measure for $\IFS(\cF)$ implies that every ergodic stationary measure is either equivalent or singular to $\leb$. These observations indicate a possible approach for Question \ref{q:q.c-branch}, particularly in dimension 2, where the work of \cite{Brown_Hertz_2017} provides a classification of all ergodic stationary measures. It is natural to ask the following question.
    
    \begin{question}
    Given a generic pair of nearby $\cC^s$ Anosov diffeomorphisms on $\bT^2$ ($s>1$), does there exist a distribution on the generated semigroup for which there is some ergodic stationary measure in the 
    Lebesgue class?
    \end{question}

    \end{document}